\documentclass[11pt]{article}

\makeatletter
\newfont{\footsc}{cmcsc10 at 8truept}
\newfont{\footbf}{cmbx10 at 8truept}
\newfont{\footrm}{cmr10 at 10truept}
\makeatother
\usepackage{amsfonts,amsmath,amsopn,amssymb,amsthm,bbm,latexsym,mathrsfs,pstricks,verbatim}
\let\n\noindent

\font\small=cmr8

\font\tenmsy=msbm10
\font\sevenmsy=msbm10 at 7pt
\font\fivemsy=msbm10 at 5pt
\newfam\msyfam % family 11
\textfont\msyfam=\tenmsy
\scriptfont\msyfam=\sevenmsy
\scriptscriptfont\msyfam=\fivemsy

\let\s\sigma

\let\l\left
\let\r\right

\def\la{{\lambda}}

\def\e{{\epsilon}}
\def\y{{\infty}}

\let\Rw\Rightarrow
\def\l{{\left}}
\def\r{{\right}}
\def\rw{\rightarrow}

\newtheorem{theorem}{Theorem}

\newtheorem{proposition}[theorem]{Proposition}

\theoremstyle{remark}

\let\s\sigma

\let\l\left
\let\r\right

\font\small=cmr8

\begin{document}

\vskip18pt

 \title{\vskip60pt Jagged partitions and lattice paths}

% \title{\vskip60pt Paths and jagged partitions: a
% configuration path description of graded parafermions % and the $\SM(2,4\ka$)  models
   %  }
\vskip18pt

%superconformal minimal

\smallskip
\author{ P. Jacob and P. Mathieu\thanks{patrick.jacob@durham.ac.uk,
pmathieu@phy.ulaval.ca.} \\
\\
Department of Mathematical Sciences, \\University of Durham, Durham, DH1 3LE, UK\\
and\\
D\'epartement de physique, de g\'enie physique et d'optique,\\
Universit\'e Laval,
Qu\'ebec, Canada, G1K 7P4.
}% \foot{Corresponding author (PM): tel:
% 418-656-3416; fax: 418-656-2040 }

\vskip .2in
\bigskip
%\bigskip
%\bigskip
%\bigskip
\date{May 2006}

\maketitle

% \begin{abstract}

\vskip0.3cm
\centerline{{\bf ABSTRACT}}
\vskip18pt

A lattice-path description of $K$-restricted jagged partitions is presented. The corresponding  lattice paths can have peaks only at even $x$ coordinate and  the maximal value of the height cannot be larger than $K-1$. Its weight is twice that of the corresponding jagged  partitions. The equivalence is demonstrated at the level of generating functions. A  bijection is given between  $K$-restricted jagged partitions and  partitions restricted by  the following frequencies conditions: $f_{2j-1}$ is even and  $f_j+f_{j+1}\leq K-1$, where $f_j$ is the number of occurrences of the part $j$ in the partition.
Bijections are given  between  paths and these restricted partitions and between  paths  and partitions  with  successive ranks in a prescribed interval.

\newpage

%==============================================================================
\let\Rw\Rightarrow
\let\rw\rightarrow
\let\l\left
\let\r\right
\let\s\sigma
\let\ka\kappa
\let\de\delta
%=================================================================

\section{Introduction}

A jagged partition of length $m$ is a sequence of non-negative integers $(n_1,\cdots , n_m)$ such that
\begin{equation}\label{jag}
 n_j\geq n_{j+1}-1\;,\qquad  \qquad  n_j\geq n_{j+2}\;, \qquad\qquad   n_m\geq 1\;.\
\end{equation}
It differs from an ordinary  partition in that there is a possible increase of one between adjacent parts $n_j$ and $n_{j+1} $. For instance $(3,4,3,2,1,2,1,0,1)$ is a jagged partition of 17.
A $K$-restricted jagged partition is a jagged partition that is further subject to  the  following restrictions:
\begin{equation}\label{rjag}
n_j \geq  n_{j+K-1} +1 \qquad{\rm or} \qquad
 n_j = n_{j+1}-1 =  n_{j+K-2}+1= n_{j+K-1} \;,
\end{equation}
for all values of $j\leq m-K+1$, with $K>1$.

Jagged partitions have been introduced in \cite{JM} in the framework of a problem in conformal field theory. In that context, the restriction reflects a generalized exclusion principle. They have been  further studied in \cite{FJM.R}. The  generating functions for restricted jagged partitions has been presented in  \cite{FJM.E} (a  special case having   been treated previously in  \cite{BFJM}). This construction relies on a recurrence method very similar to the one used by Andrews \cite{Andrr, Andr} to obtain the multiple sum expression enumerating the partitions whose parts are subject to the conditions $\la_i\geq \la_{i+k-1}+2$.

 In the conclusion of \cite{FJM.R}, the problem  of the lattice-path description of restricted jagged  partitions has been raised.
 %(it  was actually suggested by the referee).
 In the present article, we provide such a description. The main result is the following:

 \begin{theorem}\label{J=P} Let  $K,\,i,\,\ka,\,\e$ be  positive integers such that  $K>2$, $K=2\ka-\e, $ $\e=0,1$, and  $1\leq i\leq \ka$. Then the following two sets of  numbers are equal:

\n $J_{K,i}(n,m)$: the number of sequences of $m$ non-negative integers $(n_1,\cdots , n_m)$ satisfying (\ref{jag})-(\ref{rjag}) with $n= \sum_jn_j$ and  containing at most $i-1$ pairs of 01;

\n $P_{K,i}(2n,m)$:  the number of lattice paths  of weight $2n$,  charge $m$, starting at $(0,2\ka-2i)$, with peaks of height not larger than $K-1$ and all peaks occurring at  even $x$ positions.
%Then, we have $J_{K,i}(n,m) = P_{K,i}(2n,m)$.
\end{theorem}

This result relies on a  new form of the generating function of $J_{K,i}(n,m)$ given in Sect. 2. The restricted paths enumerated by $P_{K,i}(2n,m)$ are introduced in Sect. 3. The concepts of weight and charge for paths are defined there. The equality  $J_{K,i}(n,m) = P_{K,i}(2n,m)$ is demonstrated in Sect. 4 via the equivalence of their generating function.

This result allows us to connect the problem of enumerating
restricted jagged partitions to  other partition problems. These results  (adding a little piece from \cite{FJM.E}) are summarized in the following theorem.

\begin{theorem}\label{burge} Using the numbers  introduced in Theorem \ref{J=P},  with $J_{K,i}(n) =\sum_mJ_{K,i}(n,m)$, and

\n $E_{K,i}(2n,m)$:  the number of partitions $(p_1,\cdots , p_m)$  of $2n$ into $m$ parts such that very odd part occurs an even  number of times, satisfying $p_i\geq p_{i+K-1}+2$ and  containing at  most $2i-2$ copies of 1.

 \n  $R_{K,i}(2n)$:  the number of partitions   of $2n$ with all successive ranks being odd and lying in the range $[3-2i,\, 2K-1-2i]$.

  \n $O_{K,i}(n)$: the number of overpartitions
  of $n$   where non-overlined parts are not equal to
 $0,\pm i$ {\rm mod} $ K+1$ if
 $1\leq 2i< K+1$ or  where no part  is equal to 0 {\rm mod}  $ \kappa$ if $\e=1$ and
 $i=\kappa$.

\n Then, we have
$J_{K,i}(n,m) = E_{K,i}(2n,m) $ and $J_{K,i}(n)= R_{K,i}(2n)= O_{K,i}(n)$.
\end{theorem}

The equality $J_{K,i}(n,m) = E_{K,i}(2n,m)$   is proved in Sect. 5.   Bijective proofs of the equalities  $P_{K,i}(2n)=E_{K,i}(2n)$ and $P_{K,i}(2n) = R_{K,i}(2n)$ are presented respectively in Sects 6 and 7. Finally, the Rogers-Ramanujan-type identity underlying the proof of  the equality $J_{K,i}(n)   = O_{K,i}(n)$, obtained previously in  \cite{FJM.E} (Sect. 6), is presented in a different form in Sect. 8. Let us recall that an overpartition is a partition where  the final occurrence of a part can be overlined \cite{CL}.

 Recently, a lattice-path representation of overpartitions has been presented in \cite{CM}. Some remarks on this work and its relation to the present one are collected in the final section.

The results of this article rely to a large extend on those presented in \cite{BreL}. In particular, the constructions of Sects 4, 7 and part of 6 are variations and/or deformations of the original Bressoud's arguments. It is fair to recall that these in turn are clever reformulations of correspondences  originally obtained by  Burge \cite{Bu, Bur}, reframed in terms of binary words in \cite{AnB}. Actually, the equality $E_{K,i}(2n) = R_{K,i}(2n)$  is contained in Theorem 3 in \cite{Bur} (with a slight correction).  In this perspective, the above Theorem {\ref{burge} provides a new link, namely to jagged partitions.

%=================================================================

\section{New generating function for restricted  jagged partitions}

Let us first recall the expression for the generating function for $K$-restricted jagged  partitions \cite{FJM.E} (Theorems 1 and 7). Its expression uses
the notation
\begin{equation*}
(a)_n=(a;q)_n= \prod_{i=0}^{n-1} (1-aq^i)\; .
\end{equation*}

\begin{proposition}\label{J}  With $J_{K,i}(n,m)$ defined in Theorem 1 and
\begin{equation*}
J_{K,i}(z;q)= \sum_{n,m} J_{K,i}(n,m) q^n z^m\,,
\end{equation*}
we have
\begin{equation*}\label{grande}
J_{K,i}(z;q)= \sum_{m_1,\cdots,m_{\ka-1}=0}^\y    \frac{ (-zq^{1+\e m_{\ka-1}})_\y\,
q^{N_1^2+\cdots+ N_{\ka-1}^2+L_{i}}\; z^{2N} }{ (q)_{m_1}\cdots (q)_{m_{\ka-1}}
}\;,
\end{equation*}
with
\begin{equation*}\label{defNL}
 N_j= m_j+\cdots
+m_{\ka-1}\,,\qquad L_j=N_j+\cdots N_{\ka-1}\,,\qquad N=L_1\; .
\end{equation*}
\end{proposition}

An alternative expression for the generating function of $K$-restricted  jagged partitions is as follows.

\begin{proposition} \label{GFjp}
Introduce
\begin{equation*}\label{defFabb}
G_{K,i}(z; q)= \sum_{m_{1}, m_2,\cdots,m_{K-1}=0}^\y  \frac{q^{\frac12({\tilde N}_{1}^2+{\tilde N}_2^2+\cdots+ {\tilde N}_{K-1}^2 +{\tilde M})+{\tilde L}_{2i}} z^{{\tilde N}}
}{ (q)_{m_{1}}(q)_{m_2}\cdots (q)_{m_{K-1}} } \;,
\end{equation*}
where
%\begin{equation*}
\begin{align*}
 {\tilde N}_{j}&= m_j+m_{j+1}+\cdots + m_{K-1} &
{\tilde L}_{2i} &= m_{2i}+m_{2i+1}+2(m_{2i+2}+m_{2i+3})+\cdots  \cr
{\tilde M} &= m_{1}+m_{3}+\cdots +m_{K-1-\e}  &
{\tilde N}&= m_{1}+2m_2+\cdots +(K-1)m_{K-1} \;.
\end{align*}
%\end{equation*}
Then we have $
 G_{K,i}(z; q)= J_{K,i}(z; q) $,
 with $J_{K,i}(z; q)$ defined in Proposition \ref{J}.
 \end{proposition}

\begin{proof} We confine ourself to a sketch of an analytic  proof. The trick is to sum up exactly the odd modes of $G_{K,i}(z; q)$ one by one, starting form the largest one down to $m_3$. Each summation is performed by  the $q$-binomial theorem. Then $m_1$ is summed  with the Euler relation. Finally, redefining the modes $m_{2j}$ as $m_j$, we precisely recover $J_{K,i}(z;q)$.

 Let us illustrate the  argument for $K=5$, for which  $G_{5,i}(z;q) $ reads
\begin{equation*}
G_{5,i}(z;q) = \sum_{m_i=0 }^\infty    \frac{ q^{\Delta}  \, z^{ m_1+2m_2+3m_3+4m_4} }{
(q)_{m_1}(q)_{m_2}(q)_{m_3} (q)_{m_4} }\;,
\end{equation*}
with
\begin{equation*}
\begin{split}
\Delta=  &\tfrac12(m_1+m_2+m_3+m_4)^2+\tfrac12(m_2+m_3+m_4)^2+\tfrac12(m_3+m_4)^2+\tfrac12m_4^2\cr
& +\tfrac12(m_1+m_3)+(2-i)(m_2+m_3)+(3-i)m_4\; .
\end{split}
\end{equation*}
We first  replace $m_2$ by $m_2-m_3$ (with the convention that $1/(q)_n=0$ if  $n<0$) and
 use the  $q$-binomial theorem (\cite{Andr}, eq (3.3.6))
%  \begin{equation}\sum_{j=0}^n {(q)_n\over
% (q)_j (q)_{n-j}} q^{\frac12j(j+1)} x^j = (-xq)_n \;,\
%  \end{equation}
to  perform the summation over $m_3$:
 \begin{equation*}
\frac1{(q)_{m_2} } \sum_{m_3=0}^{m_2} \frac{ (q)_{m_2} }{ (q)_{m_3} (q)_{m_2-m_3}  }    q^{\frac12m_3(m_3+1)} (zq^{m_4})^{m_3} = \frac{(-zq^{1+m_4})_{m_2} }{ (q)_{m_2} }
\end{equation*}
We have thus obtained at this point the form
 \begin{equation*}
 G_{5,i}(z;q)=
 \sum_{m_1,m_2,m_4=0 }^\infty  \frac{ (-zq^{1+m_4})_{m_2}  q^{\Delta'} z^{m_1+2m_2+4m_4}}{
 (q)_{m_1}(q)_{m_2}(q)_{m_4} }
 \end{equation*}
 with
  \begin{equation*}
  \Delta'= \frac12(m_1+m_2+m_4)^2+\frac12(m_2+m_4)^2+\frac12m_4^2+\frac12m_1+ (2-i)m_2+(3-i)m_4
  \end{equation*}
% = \sum_{m_1,m_2,m_3=0 }^\infty { q^{\frac12m_1(m_1+2m_2+1)
% +m_2^2+\frac12(m_3(m_3+1) + (2-i)m_2} z^{m_1+2m_2+m_3}\over
% (q)_{m_1}(q)_{m_2-m_3}(q)_{m_3}}\cr &
%  & =\sum_{m_1,m_2=0 }^\infty {q^{\frac12m_1(m_1+2m_2+1)+m_2^2+(2-i)m_2} z^{m_1+2m_2} \over
% (q)_{m_1}(q)_{m_2}} \sum_{m_3=0}^{m_2} {(q)_{m_2}\over
% (q)_{m_2-m_3}(q)_{m_3}} q^{\frac12m_3(m_3+1) } z^{m_3}\;.\cr}\eq
% $$
Next, we use the Euler relation  (\cite{Andr}, eq (2.2.6))
% $\sum{q^{\frac12n(n-1)} x^n (q)_n^{-1}} = (-x)_\infty$
  \begin{equation*}
 \sum_{n=0 }^\infty{q^{\frac12n(n-1)} x^n\over (q)_n} = (-x)_\infty\;,
  \end{equation*}
 to sum over  $m_1$. This results into
 \begin{equation*}
 G_{5,i}(z;q)=
\sum_{m_2,m_4=0}^\infty  \frac{ (-zq^{1+m_4})_\infty q^{(m_2+m_4)^2+m_4^2+(2-i)m_2+(3-i)m_4}z^{m_2+2m_4}  }{
(q)_{m_2}(q)_{m_4}}\;.\end{equation*}
With $m_{2j}$ redefined as $m_j$, this is precisely $J_{5,i}(z;q)$ as given in Proposition (\ref{J}).
 \end{proof}

\n {\it Remark}: The idea of this proof is due to Warnaar \cite{Wa}. For $\e=0$, this argument has been presented in Appendix  A of \cite{JM.A} where this generating function has  already been displayed.
In its full generality, this can be proved from  a extension of Lemma A.1 of \cite{War}. A combinatorial proof of $G_{K,i}(1; q^2)= J_{K,i}(1; q^2) $ in the case $\e=1$ is given in \cite{Bur} (cf. the last equation of p. 204). Finally, an alternative combinatorial proof of Proposition \ref{GFjp} pertaining to all values of $K>2$ follows from combining the results of Sects 4-7 below.

 % \n Preuve directe si on a la bijection? Voir Burge et Foada et al.

%=================================================================
\section{Restricted  lattice paths}

We now  introduce our lattice paths as natural generalizations of the Bressoud's ones \cite{BreL}.  Paths are defined in the first quadrant of   an integer square lattice by the following rules:

\n 1- A path starts at an even non-negative  integer  position on the $y$ axis.
% , i.e., at $(0,j)$ with $j$ a non-negative integer.

\n 2- The possible moves are
either  from $(i,j)$ to $(i+1,j+1)$ or  from $(i,j)$ to $(i+1,j-1)$ if $j\geq 1$  or from $(i,0)$ to $(i+1,0)$. (The horizontal move is thus  allowed only on the $x$ axis.)

\n 3- A non-empty path always terminates on the $x$ axis, i.e., with the move $(i,1)$ to $(i+1,0)$.

\n 4- The peaks can only be at even integer positions.

\n 5- The height (i.e., $y$ coordinate) of the peaks cannot be larger that $K-1$.

\n The essential difference between these paths and those of Bressoud is that here the initial vertical position as well as the peak positions are forced to be even.

We next need to define the concepts of {\it weight}, {\it  charge} and {\it relative height}. The weight of a path is simply the sum of the $x$ coordinate of all the peaks.  The relative height of a peak with $(x,y)$ coordinate $(i,j)$ is the largest integer $h$ such that we can find two vertices $(i',j-h)$ and $(i'',j-h)$ on the path  with $i'<i<i''$ and such that between these two vertices there are no peaks of height larger than $j$ and every peak of height equal to $j$ have weight larger than $i$ \cite{BP}.
The charge of a path is the sum of all its relative heights.

Denote by the pair $(x_j;h_j)$ the position and relative height of the $j$-th  peak counted from the left.  A path is fully determined by the data $\{(x_j;h_j ) \}$ and the initial position $a$. An example is presented in Figure 1.

% , which specify the position of the peaks $x_j$ and their relative height $h_j$.

\n {\it Remark:} Similar but slightly different paths (horizontal moves not being  allowed) have been  introduced in \cite{Ole}. There the charge stands for  the relative height. Our notion of charge is thus different from this one as it refers to the sum of all the individual charges of a path.

%\n {\bf Figure 1: un exemple de chemin avec la determination des hauteurs relatives}

\begin{figure}[ht]
\caption{{\footnotesize A path of weight $56$ starting at $a=2$ and with peaks at
positions $2$, $6$, $10$, $14$ and $24$.  The relative heights are
respectively $1$, $1$, $4$, $2$ and $5$.  This is a $K$-restricted
path for all $K \geq 6$.}} \label{figure1}
\begin{center}
\begin{pspicture}(0,0)(15.5,4)
%axis
\psline{->}(0.5,0.5)(0.5,3.5) \psline{->}(0.5,0.5)(15.5,0.5)
%units
\psline{-}(0.5,0.5)(0.5,0.6) \psline{-}(1.0,0.5)(1.0,0.6)
\psline{-}(1.5,0.5)(1.5,0.6) \psline{-}(2.0,0.5)(2.0,0.6)
\psline{-}(2.5,0.5)(2.5,0.6) \psline{-}(3.0,0.5)(3.0,0.6)
\psline{-}(3.5,0.5)(3.5,0.6) \psline{-}(4.0,0.5)(4.0,0.6)
\psline{-}(4.5,0.5)(4.5,0.6) \psline{-}(5.0,0.5)(5.0,0.6)
\psline{-}(5.5,0.5)(5.5,0.6) \psline{-}(6.0,0.5)(6.0,0.6)
\psline{-}(6.5,0.5)(6.5,0.6) \psline{-}(7.0,0.5)(7.0,0.6)
\psline{-}(7.5,0.5)(7.5,0.6) \psline{-}(8.0,0.5)(8.0,0.6)
\psline{-}(8.5,0.5)(8.5,0.6) \psline{-}(9.0,0.5)(9.0,0.6)
\psline{-}(9.5,0.5)(9.5,0.6) \psline{-}(10.0,0.5)(10.0,0.6)
\psline{-}(10.5,0.5)(10.5,0.6) \psline{-}(11.0,0.5)(11.0,0.6)
\psline{-}(11.5,0.5)(11.5,0.6) \psline{-}(12.0,0.5)(12.0,0.6)
\psline{-}(12.5,0.5)(12.5,0.6) \psline{-}(13.0,0.5)(13.0,0.6)
\psline{-}(13.5,0.5)(13.5,0.6) \psline{-}(14.0,0.5)(14.0,0.6)
\psline{-}(14.0,0.5)(14.0,0.6)
\psline{-}(14.5,0.5)(14.5,0.6)\rput(1.5,0.25){{\small $2$}}
\rput(3.5,0.25){{\small $6$}} \rput(5.5,0.25){{\small $10$}}
\rput(7.5,0.25){{\small $14$}} \rput(12.5,0.25){{\small $24$}}
 \psline{-}(0.5,1.0)(0.6,1.0)
\psline{-}(0.5,1.5)(0.6,1.5) \psline{-}(0.5,2.0)(0.6,2.0)
\psline{-}(0.5,2.5)(0.6,2.5) \psline{-}(0.5,3.0)(0.6,3.0)
\rput(-0.08,1.5){{\small $a=2$}} \rput(0.25,2.5){{\small $4$}}
%graphic
\psline{-}(0.5,1.5)(1.0,1.0) \psline{-}(1.0,1.0)(1.5,1.5)
\psline{-}(1.5,1.5)(2.0,1.0) \psline{-}(2.0,1.0)(2.5,0.5)
\psline{-}(2.5,0.5)(3.0,1.0) \psline{-}(3.0,1.0)(3.5,1.5)
\psline{-}(3.5,1.5)(4.0,1.0) \psline{-}(4.0,1.0)(4.5,1.5)
\psline{-}(4.5,1.5)(5.0,2.0) \psline{-}(5.0,2.0)(5.5,2.5)
\psline{-}(5.5,2.5)(6.0,2.0) \psline{-}(6.0,2.0)(6.5,1.5)
\psline{-}(6.5,1.5)(7.0,2.0) \psline{-}(7.0,2.0)(7.5,2.5)
\psline{-}(7.5,2.5)(8.0,2.0) \psline{-}(8.0,2.0)(8.5,1.5)
\psline{-}(8.5,1.5)(9.0,1.0) \psline{-}(9.0,1.0)(9.5,0.5)
\psline{-}(9.5,0.5)(10.0,0.5)
 \psline{-}(10.0,0.5)(10.5,1.0)\psline{-}(10.5,1.0)(11.0,1.5)
\psline{-}(11.0,1.5)(11.5,2.0) \psline{-}(11.5,2.0)(12.0,2.5)
\psline{-}(12.0,2.5)(12.5,3.0) \psline{-}(12.5,3.0)(13.0,2.5)
\psline{-}(13.0,2.5)(13.5,2.0) \psline{-}(13.5,2.0)(14.0,1.5)
\psline{-}(14.0,1.5)(14.5,1.0) \psline{-}(14.5,1.0)(15.0,0.5)

\end{pspicture}
\end{center}
\end{figure}
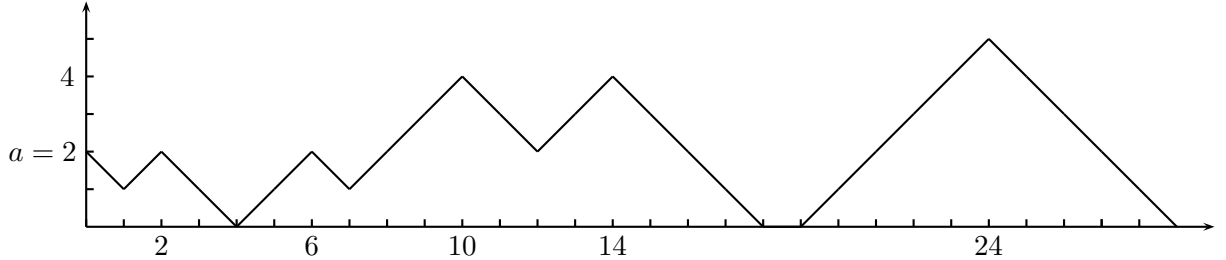

%\rput(1.0,0.25){{\small $-L=-3$}} \rput(2.0,2.0){{\small$\ell=3$}} \rput(10.5,0.25){{\small $M=16$}}
%=================================================================

\section{Generating function for restricted paths}

This section is devoted to obtaining the generating function  for  our lattice paths. More precisely, we verify in the following proposition that $G_{K,i}(z;q^2)$ is the generating function  for those paths enumerated by $P_{K,i}(2n,m)$. Propositions \ref{GFjp}   and \ref{GFpa} prove Theorem \ref{J=P}. Note that introducing the  charge allows us to extend the lattice-path  generating function  to a two-variable function.

\begin{proposition}\label{GFpa}  The generating function for $K$-restricted paths starting at the vertical position $(0,2\ka-2i)$, with  $1\leq i\leq \ka$,   is given by  $G_{K,i}(z;q^2)$ defined in Proposition \ref{GFjp}. In this multiple-sum expression, $m_j$ is the number of peaks with relative height $j$ and the power of $z$ is the charge of the path.
\end{proposition}

\begin{proof}
We use an inductive argument on $K$ which is a simple adaptation of the Bressoud's proof (cf. \cite{BreL} Sect. 2). However, since the inductive step is $K\rw K+2$, we first need
 to verify the result  for $K=2$ and $3$.

  Let us then consider $K=2$. The maximal height is $1$. In that case all paths necessarily  start at $(0,0)$. For  a path containing $m$ adjacent peaks of height $1$, the minimal-weight configuration is obtained by an initial  horizontal move of length $1$ followed by  a sequence of $m$ peaks in contact, with $x$ coordinates $2,4,\dots, 2m$. The weight of this configuration is $m(m+1)$.  Next, the peaks can be moved along the $x$ axis, starting for the right-most one, and proceeding successively from right to left. Denote the respective displacements by the sequence  of even integers $(\ell_1,\ell_2,\cdots, \ell_m)$ with $\ell_i\geq \ell_{i+1} $ and $\ell_i\geq 0$, where $\ell_1$ is the displacement of the rightmost peak, $\ell_2$ that of the next one to its left, etc. Such sequences are partitions containing  at most $m$ even parts. Their generating function is $1/(q^2;q^2)_m$.  Therefore, the generating function for all paths starting at (0,0) and containing $m$ peaks of height $1$ (the charge being thus $m$) is
\begin{equation*}\frac{q^{m(m+1)} z^m }{ (q^2;q^2)_m} \;.
\end{equation*}
By summing over $m$, one recovers $G_{2,1}(z;q^2)$.

Consider next  $K=3$, so that $\ka=2$ and $i=1,\,2$. The maximal height being now $2$, the paths can have peaks of relative heights $1$ or 2. Let us start with the case $i=2$, that is,  paths starting at the origin.  One first  needs to find the minimal-weight configuration, starting at (0,0), that contains $m_{1}$ peaks of relative height $1$ and  $m_2$ peaks of relative height 2. This can be determined by comparing the lowest  weight associated to all possible ordering of the peaks.
In this way, we find that the minimal-weight configuration with height-content $m_{1}$ and $m_2$ is obtained with a  peak of height 2 at position 2, followed by $m_{1}$ successive peaks of  relative height $1$ (but height 2), at position $4,6,\cdots, 2m_{1}+2$,  and then $m_{2}-1 $ successive peaks of height 2, at position $2m_{1}+6,Ê\cdots, 2m_{1}+4(m_2-2)+2$. This path has total weight
\begin{equation*}
w_0^{(i=2)}= m_{1}(m_{1}+1) +2m_2^2+2m_{1}m_2\,.
\end{equation*}
The minimal-weight configuration for a path starting at the origin and with $(m_{1},m_2)=(3,4)$ is presented in Figure 2.

%\n {\bf Figure 2: configuration minimale pour un depart a (0,0) pour   $(m_{1},m_2)=(3,4)$.}

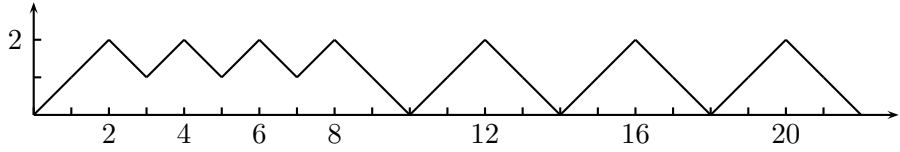
\begin{figure}[ht]
\caption{{\footnotesize Minimal weight configuration for $\kappa=2$,
$i=2$ and $(m_{1},m_2)=(3,4)$.}} \label{figure2}
\begin{center}
\begin{pspicture}(0,0)(12.0,2.5)
%axis
\psline{->}(0.5,0.5)(0.5,2.0) \psline{->}(0.5,0.5)(12.0,0.5)
%units
\psline{-}(0.5,0.5)(0.5,0.6) \psline{-}(1.0,0.5)(1.0,0.6)
\psline{-}(1.5,0.5)(1.5,0.6) \psline{-}(2.0,0.5)(2.0,0.6)
\psline{-}(2.5,0.5)(2.5,0.6) \psline{-}(3.0,0.5)(3.0,0.6)
\psline{-}(3.5,0.5)(3.5,0.6) \psline{-}(4.0,0.5)(4.0,0.6)
\psline{-}(4.5,0.5)(4.5,0.6) \psline{-}(5.0,0.5)(5.0,0.6)
\psline{-}(5.5,0.5)(5.5,0.6) \psline{-}(6.0,0.5)(6.0,0.6)
\psline{-}(6.5,0.5)(6.5,0.6) \psline{-}(7.0,0.5)(7.0,0.6)
\psline{-}(7.5,0.5)(7.5,0.6) \psline{-}(8.0,0.5)(8.0,0.6)
\psline{-}(8.5,0.5)(8.5,0.6) \psline{-}(9.0,0.5)(9.0,0.6)
\psline{-}(9.5,0.5)(9.5,0.6) \psline{-}(10.0,0.5)(10.0,0.6)
\psline{-}(10.5,0.5)(10.5,0.6) \psline{-}(11.0,0.5)(11.0,0.6)
%\psline{-}(11.5,0.5)(11.5,0.6) \psline{-}(12.0,0.5)(12.0,0.6)
%\psline{-}(12.5,0.5)(12.5,0.6) \psline{-}(13.0,0.5)(13.0,0.6)
%\psline{-}(13.5,0.5)(13.5,0.6) \psline{-}(14.0,0.5)(14.0,0.6)
%\psline{-}(14.0,0.5)(14.0,0.6)
%\psline{-}(14.5,0.5)(14.5,0.6)
\rput(1.5,0.25){{\small $2$}} \rput(3.5,0.25){{\small $6$}}
\rput(6.5,0.25){{\small $12$}} \rput(8.5,0.25){{\small $16$}}
\rput(2.5,0.25){{\small $4$}}\rput(4.5,0.25){{\small
$8$}}\rput(10.5,0.25){{\small $20$}}
 \psline{-}(0.5,1.0)(0.6,1.0)
\psline{-}(0.5,1.5)(0.6,1.5) %\psline{-}(0.5,2.0)(0.6,2.0)
%\psline{-}(0.5,2.5)(0.6,2.5) \psline{-}(0.5,3.0)(0.6,3.0)
\rput(0.25,1.5){{\small $2$}} %\rput(0.25,2.5){{\small $2 $}}
%graphic
\psline{-}(0.5,0.5)(1.0,1.0) \psline{-}(1.0,1.0)(1.5,1.5)
\psline{-}(1.5,1.5)(2.0,1.0) \psline{-}(2.0,1.0)(2.5,1.5)
\psline{-}(2.5,1.5)(3.0,1.0) \psline{-}(3.0,1.0)(3.5,1.5)
\psline{-}(3.5,1.5)(4.0,1.0) \psline{-}(4.0,1.0)(4.5,1.5)
\psline{-}(4.5,1.5)(5.0,1.0) \psline{-}(5.0,1.0)(5.5,0.5)
\psline{-}(5.5,0.5)(6.0,1.0) \psline{-}(6.0,1.0)(6.5,1.5)
\psline{-}(6.5,1.5)(7.0,1.0) \psline{-}(7.0,1.0)(7.5,0.5)
\psline{-}(7.5,0.5)(8.0,1.0) \psline{-}(8.0,1.0)(8.5,1.5)
\psline{-}(8.5,1.5)(9.0,1.0) \psline{-}(9.0,1.0)(9.5,0.5)
\psline{-}(9.5,0.5)(10.0,1.0)
 \psline{-}(10.0,1.0)(10.5,1.5)\psline{-}(10.5,1.5)(11.0,1.0)
\psline{-}(11.0,1.0)(11.5,0.5)% \psline{-}(11.5,2.0)(12.0,2.5)
%\psline{-}(12.0,2.5)(12.5,3.0) \psline{-}(12.5,3.0)(13.0,2.5)
%\psline{-}(13.0,2.5)(13.5,2.0) \psline{-}(13.5,2.0)(14.0,1.5)
%\psline{-}(14.0,1.5)(14.5,1.0) \psline{-}(14.5,1.0)(15.0,0.5)

\end{pspicture}
\end{center}
\end{figure}

The peaks can then be displaced by moves to be described in more detail below. It suffices at this point to note that  the displacement of the  peaks  of relative height $1$ and 2 are determined by two partitions of even integers with at most $m_{1}$ and $m_{2}$ even parts respectively.  The total charge is $m_1+2m_2$. The generating function for these paths is thus
\begin{equation*}
\frac{ q^{m_{1}(m_{1}+1) +2m_2^2+2m_{1}m_2} \;  z^{m_1+2m_2}}{(q^2;q^2)_{m_{1} }(q^2;q^2)_{m_{2}}}\;.
\end{equation*}
The summation over $ m_{1}$ and $ m_2$ yields $G_{3,2}(z;q^2)$.

Consider next $i=1$, that is,  paths  starting at $(0,2)$. The minimal-weight configuration is slightly modified in that case; it is given by  $m_{1}$ successive peaks of  relative height $1$ (and height 2), at position $2,\, 4,\cdots, 2m_{1}$ (after which   the path joins the horizontal axis), followed by $m_{2} $ successive peaks of height 2, at position $2m_{1}+4Ê \cdots, \,2 m_{1}+4m_2$. This is illustrated in Figure 3. The corresponding weight is
\begin{equation*}
w_{0}^{(i=1)}= m_{1}(m_{1}+1) +2m_2^2+2m_{1}m_2 +2m_2= w_0^{(i=2)}+2m_2\,.
\end{equation*}
The generating function is constructed as before and it reads
\begin{equation*} \frac{  q^{m_{1}(m_{1}+1) +2m_2^2+2m_{1}m_2 +2m_2 } \; z^{m_1+2m_2}} {(q^2;q^2)_{m_{1} }(q^2;q^2)_{m_{2}}}\; .
\end{equation*} The summation over $m_{1} $ and $m_2$ reproduces $G_{3,1}(z;q^2)$.

 %

 %\footnote{This way of probing the linear terms will be made more systematic below.}

%\n {\bf Figure 3: configuration minimale pour un depart a (0,1) pour $(m_{1/2}, m_1)=(3,4)$.}

\begin{figure}[ht]
\caption{{\footnotesize Minimal weight configuration for $\kappa=2$,
$i=1$ and $(m_{1},m_2)=(3,4)$.}} \label{figure3}
\begin{center}
\begin{pspicture}(0,0)(13.0,2.5)
%axis
\psline{->}(0.5,0.5)(0.5,2.0) \psline{->}(0.5,0.5)(13.0,0.5)
%units
\psline{-}(0.5,0.5)(0.5,0.6) \psline{-}(1.0,0.5)(1.0,0.6)
\psline{-}(1.5,0.5)(1.5,0.6) \psline{-}(2.0,0.5)(2.0,0.6)
\psline{-}(2.5,0.5)(2.5,0.6) \psline{-}(3.0,0.5)(3.0,0.6)
\psline{-}(3.5,0.5)(3.5,0.6) \psline{-}(4.0,0.5)(4.0,0.6)
\psline{-}(4.5,0.5)(4.5,0.6) \psline{-}(5.0,0.5)(5.0,0.6)
\psline{-}(5.5,0.5)(5.5,0.6) \psline{-}(6.0,0.5)(6.0,0.6)
\psline{-}(6.5,0.5)(6.5,0.6) \psline{-}(7.0,0.5)(7.0,0.6)
\psline{-}(7.5,0.5)(7.5,0.6) \psline{-}(8.0,0.5)(8.0,0.6)
\psline{-}(8.5,0.5)(8.5,0.6) \psline{-}(9.0,0.5)(9.0,0.6)
\psline{-}(9.5,0.5)(9.5,0.6) \psline{-}(10.0,0.5)(10.0,0.6)
\psline{-}(10.5,0.5)(10.5,0.6) \psline{-}(11.0,0.5)(11.0,0.6)
\psline{-}(11.5,0.5)(11.5,0.6) \psline{-}(12.0,0.5)(12.0,0.6)
%\psline{-}(12.5,0.5)(12.5,0.6) \psline{-}(13.0,0.5)(13.0,0.6)
%\psline{-}(13.5,0.5)(13.5,0.6) \psline{-}(14.0,0.5)(14.0,0.6)
%\psline{-}(14.0,0.5)(14.0,0.6)
%\psline{-}(14.5,0.5)(14.5,0.6)
\rput(1.5,0.25){{\small $2$}} \rput(3.5,0.25){{\small $6$}}
\rput(5.5,0.25){{\small $10$}} \rput(7.5,0.25){{\small $14$}}
\rput(2.5,0.25){{\small $4$}}\rput(9.5,0.25){{\small
$18$}}\rput(11.5,0.25){{\small $22$}}
 \psline{-}(0.5,1.0)(0.6,1.0)
\psline{-}(0.5,1.5)(0.6,1.5) %\psline{-}(0.5,2.0)(0.6,2.0)
%\psline{-}(0.5,2.5)(0.6,2.5) %\psline{-}(0.5,3.0)(0.6,3.0)
\rput(0.25,1.5){{\small $2$}} %\rput(0.25,2.5){{\small $2 $}}
%graphic
\psline{-}(0.5,1.5)(1.0,1.0) \psline{-}(1.0,1.0)(1.5,1.5)
\psline{-}(1.5,1.5)(2.0,1.0) \psline{-}(2.0,1.0)(2.5,1.5)
\psline{-}(2.5,1.5)(3.0,1.0) \psline{-}(3.0,1.0)(3.5,1.5)
\psline{-}(3.5,1.5)(4.0,1.0) \psline{-}(4.0,1.0)(4.5,0.5)
\psline{-}(4.5,0.5)(5.0,1.0) \psline{-}(5.0,1.0)(5.5,1.5)
\psline{-}(5.5,1.5)(6.0,1.0) \psline{-}(6.0,1.0)(6.5,0.5)
\psline{-}(6.5,0.5)(7.0,1.0) \psline{-}(7.0,1.0)(7.5,1.5)
\psline{-}(7.5,1.5)(8.0,1.0) \psline{-}(8.0,1.0)(8.5,0.5)
\psline{-}(8.5,0.5)(9.0,1.0) \psline{-}(9.0,1.0)(9.5,1.5)
\psline{-}(9.5,1.5)(10.0,1.0)
 \psline{-}(10.0,1.0)(10.5,0.5)\psline{-}(10.5,0.5)(11.0,1.0)
\psline{-}(11.0,1.0)(11.5,1.5) \psline{-}(11.5,1.5)(12.0,1.0)
\psline{-}(12.0,1.0)(12.5,0.5)% \psline{-}(12.5,3.0)(13.0,2.5)
%\psline{-}(13.0,2.5)(13.5,2.0) \psline{-}(13.5,2.0)(14.0,1.5)
%\psline{-}(14.0,1.5)(14.5,1.0) \psline{-}(14.5,1.0)(15.0,0.5)

\end{pspicture}
\end{center}
\end{figure}

We are now in position to work out the inductive step. Let us suppose that the result is true for $K$ and prove it for $K+2$. The hypothesis is thus that the generating function for lattice paths that start at $(0,2\ka+2-2i)$ with height $\leq K-1$ and having  $n_j$  peaks of relative height at least $ j-1$ is given by:
 \begin{equation}\label{hypo}
\frac{q^{({ n}_{3}^2+{ n}_4^2+\cdots+ { n}_{K+1}^2) +(n_{3}-n_4+n_{5}-n_6\cdots
+n_{K+1-\e}-\e n_{K+2-\e})+2(n_{2i}+n_{2i+2}+\cdots )}\; z^{n_3+n_4+\cdots+n_{K+1}} }{  (q^2;q^2)_{n_{3}-n_4}(q^2;q^2)_{n_4-n_{5}}\cdots (q^2;q^2)_{n_{K+1}} }\; ,
 \end{equation}
with $m_j=n_j-n_{j+1}$. Note that  we have replaced $m_j\rw m_{j+2}$ to transfer the to-be-added modes to the lowest values (1 and 2) of $j$.

The first step amounts to perform a `volcanic uplift', meaning  breaking each peak, opening the path by four  units (toward the right) and  closing the path by inserting a new peak of height 2 -- to preserve the parity of the position of the peak. (This operation is the reason for which the induction amounts to  increase $K$ by two units). This modifies the position of the first peak by 2, that of the second by 6, etc. The weight of the path is thus changed by
 \begin{equation*} w_1= 2n_{3}^2\; ,
  \end{equation*} since $n_{3}$ is the total number of peaks. This uplift has modified the upper bound of the height from $K-1$ to $K+1$.

In the second step, we insert $m_{1}$ peaks of height $1$ and  $m_2$ peaks of height 2 at the beginning of the path. This is done as for the minimal configuration of $K=3$ described above, pertaining to the case $i=2$ (which is the right choice because the original position on the vertical axis has already the desired one, namely $2\ka+2-2i$).  From left to right we thus insert one  peak of height 2 followed by  $m_{1}$  peaks of height $1$  and finally $m_2-1$ peaks of height 2 as closely packed as possible. The length of this path insertion (which also gives its charge) is $m_{1}+2m_2$. It produces a shift of the position of the  other peaks which  results into the following weight increase:
\begin{equation*} w_2= 2(m_{1}+2m_2) n_{3}\;.
  \end{equation*}
  Finally, one must take into account the weight of the added peaks, which  is
   \begin{equation*} w_3= m_{1}(m_{1}+1) +2m_2^2+2m_{1}m_2\;.
  \end{equation*}
  The total weight is then changed by
  \begin{equation*} w_1+w_2+w_3= n_{1}^2+n_2^2+ n_{1}-n_2\;.
  \end{equation*}
 The displacement of the  peaks  of relative height $1$ and 2 are again determined by two partitions with at most $m_{1}$ and $m_{2}$ even parts respectively. We start the displacement of the $m_2-1$  rightmost peaks of relative height 2 from right to left, and these moves are specified by the rules indicated in the Figure 4. If at a certain point, two or more peaks of the same height are in contact, we move the rightmost one in order to complete the required displacement. An example is given in Figure 5. Then the final peak of relative height 2 is displaced by a number of units. Its move is illustrated in Figure 6. Finally, the peaks  of relative height  $1$ are moved one by one (again from right to left ) and  with the rules displayed in  Figure 7.

 %\n{\bf Figures 4, 5,6,7: description des deplacements}

\begin{figure}[ht]
\caption{{\footnotesize Possible moves for a peak of relative height 2 within a path with all peaks being of relative height larger than 2.}}
\label{figure4}
\begin{center}
\begin{pspicture}(0,0)(6.0,6.0)
%axis
%\psline{->}(0.5,0.5)(0.5,2.0) \psline{->}(0.5,0.5)(13.0,0.5)
%units

\rput(0.75,4.5){{\small $i$}} %\rput(4.5,3.25){{\small $i+6$}}

%\rput(2.5,0.25){{\small $2$}}\rput(9.5,0.25){{\small$9$}}
%\rput(11.5,0.25){{\small $11$}}
 %\psline{-}(0.5,1.0)(0.6,1.0)
%\psline{-}(0.5,1.5)(0.6,1.5) %\psline{-}(0.5,2.0)(0.6,2.0)
%\psline{-}(0.5,2.5)(0.6,2.5) %\psline{-}(0.5,3.0)(0.6,3.0)
%\rput(0.25,1.5){{\small $1$}} %\rput(0.25,2.5){{\small $2 $}}
%graphic
\psline{-}(0.25,4.75)(0.5,5.0) \psline{-}(0.5,5.0)(0.75,5.25)
\psline{-}(0.75,5.25)(1.0,5.0) \psline{-}(1.0,5.0)(1.25,4.75)
\psline{-}(1.25,4.75)(1.5,4.75) \psline{-}(1.5,4.75)(1.75,4.75)
\psline{-}(1.75,4.75)(2.0,4.75) \psline{-}(2.0,4.75)(2.25,4.75)
%\psline{-}(4.5,5.0)(5.0,4.5) \psline{-}(5.0,4.5)(5.5,4.0)
%\psline{-}(5.5,4.0)(6.0,3.5) \psline{->}(6.3,4.5)(6.9,4.5)

%\psline[linestyle=dashed]{-}(8.0,4.0)(10.0,4.0)
%\psline[linestyle=dashed]{-}(3.5,1.0)(5.5,1.0)
\psline{->}(2.6,5.25)(3.2,5.25)

\rput(4.5,4.5){{\small $i+2$}}
%\rput(11.0,3.25){{\small $i+6$}}
%\rput(9.0,3.25){{\small $i+1$}} \rput(11.0,3.25){{\small $i+3$}}
\psline{-}(3.5,4.75)(3.75,4.75) \psline{-}(3.75,4.75)(4.0,4.75)
\psline{-}(4.0,4.75)(4.25,5.0) \psline{-}(4.25,5.0)(4.5,5.25)
 \psline{-}(4.5,5.25)(4.75,5.0)\psline{-}(4.75,5.0)(5.0,4.75)
\psline{-}(5.0,4.75)(5.25,4.75) \psline{-}(5.25,4.75)(5.5,4.75)
%\psline{-}(11.0,5.0)(11.5,4.5) \psline{-}(11.5,4.5)(12.0,4.0)
%\psline{-}(12.0,4.0)(12.5,3.5)% \psline{-}(13.5,2.0)(14.0,1.5)
%\psline{-}(14.0,1.5)(14.5,1.0) \psline{-}(14.5,1.0)(15.0,0.5)

%\rput(0.75,3.0){{\small $i$}}
%\rput(4.5,3.25){{\small $i+6$}}
%\rput(2.5,0.25){{\small $2$}}\rput(9.5,0.25){{\small$9$}}
%\rput(11.5,0.25){{\small $11$}}
 %\psline{-}(0.5,1.0)(0.6,1.0)
%\psline{-}(0.5,1.5)(0.6,1.5) %\psline{-}(0.5,2.0)(0.6,2.0)
%\psline{-}(0.5,2.5)(0.6,2.5) %\psline{-}(0.5,3.0)(0.6,3.0)
%\rput(0.25,1.5){{\small $1$}} %\rput(0.25,2.5){{\small $2 $}}
%graphic
\psline{-}(0.25,3.25)(0.5,3.5) \psline{-}(0.5,3.5)(0.75,3.75)
\psline{-}(0.75,3.75)(1.0,3.5) \psline{-}(1.0,3.5)(1.25,3.25)
\psline{-}(1.25,3.25)(1.5,3.5) \psline{-}(1.5,3.5)(1.75,3.75)
\psline{-}(1.75,3.75)(2.0,4.0) \psline{-}(2.0,4.0)(2.25,4.25)
%\psline{-}(4.5,5.0)(5.0,4.5) \psline{-}(5.0,4.5)(5.5,4.0)
%\psline{-}(5.5,4.0)(6.0,3.5) \psline{->}(6.3,4.5)(6.9,4.5)

%\psline[linestyle=dashed]{-}(8.0,4.0)(10.0,4.0)
%\psline[linestyle=dashed]{-}(3.5,1.0)(5.5,1.0)
\psline{->}(2.6,3.75)(3.2,3.75)

\rput(0.75,3.0){{\small $i$}} \rput(4.5,3.0){{\small $i+2$}}
%\rput(9.0,3.25){{\small $i+1$}} \rput(11.0,3.25){{\small $i+3$}}
\psline{-}(3.5,3.25)(3.75,3.5) \psline{-}(3.75,3.5)(4.0,3.75)
\psline{-}(4.0,3.75)(4.25,4.0) \psline{-}(4.25,4.0)(4.5,4.25)
 \psline{-}(4.5,4.25)(4.75,4.0)\psline{-}(4.75,4.0)(5.0,3.75)
\psline{-}(5.0,3.75)(5.25,4.0) \psline{-}(5.25,4.0)(5.5,4.25)
%\psline{-}(11.0,5.0)(11.5,4.5) \psline{-}(11.5,4.5)(12.0,4.0)
%\psline{-}(12.0,4.0)(12.5,3.5)% \psline{-}(13.5,2.0)(14.0,1.5)
%\psline{-}(14.0,1.5)(14.5,1.0) \psline{-}(14.5,1.0)(15.0,0.5)

%\psline{-}(12.5,0.5)(12.5,0.6) \psline{-}(13.0,0.5)(13.0,0.6)
%\psline{-}(13.5,0.5)(13.5,0.6) \psline{-}(14.0,0.5)(14.0,0.6)
%\psline{-}(14.0,0.5)(14.0,0.6)
%\psline{-}(14.5,0.5)(14.5,0.6)
%\rput(1.5,3.25){{\small $i$}} %\rput(4.5,3.25){{\small $i+6$}}

%\rput(2.5,0.25){{\small $2$}}\rput(9.5,0.25){{\small$9$}}
%\rput(11.5,0.25){{\small $11$}}
 %\psline{-}(0.5,1.0)(0.6,1.0)
%\psline{-}(0.5,1.5)(0.6,1.5) %\psline{-}(0.5,2.0)(0.6,2.0)
%\psline{-}(0.5,2.5)(0.6,2.5) %\psline{-}(0.5,3.0)(0.6,3.0)
%\rput(0.25,1.5){{\small $1$}} %\rput(0.25,2.5){{\small $2 $}}
%graphic
\psline{-}(0.25,1.75)(0.5,2.0) \psline{-}(0.5,2.0)(0.75,2.25)
\psline{-}(0.75,2.25)(1.0,2.0) \psline{-}(1.0,2.0)(1.25,1.75)
\psline{-}(1.25,1.75)(1.5,1.75) \psline{-}(1.5,1.75)(1.75,2.0)
\psline{-}(1.75,2.0)(2.0,2.25) \psline{-}(2.0,2.25)(2.25,2.5)
%\psline{-}(4.5,5.0)(5.0,4.5) \psline{-}(5.0,4.5)(5.5,4.0)
%\psline{-}(5.5,4.0)(6.0,3.5) \psline{->}(6.3,4.5)(6.9,4.5)

%\psline[linestyle=dashed]{-}(8.0,4.0)(10.0,4.0)
%\psline[linestyle=dashed]{-}(3.5,1.0)(5.5,1.0)
\psline{->}(2.6,2.25)(3.2,2.25) \psline{->}(2.6,0.75)(3.2,0.75)

\rput(0.75,1.5){{\small $i$}} \rput(4.5,1.5){{\small $i+2$}}
%\rput(9.0,3.25){{\small $i+1$}} \rput(11.0,3.25){{\small $i+3$}}
\psline{-}(3.5,1.75)(3.75,1.75) \psline{-}(3.75,1.75)(4.0,2.0)
\psline{-}(4.0,2.0)(4.25,2.25) \psline{-}(4.25,2.25)(4.5,2.5)
 \psline{-}(4.5,2.5)(4.75,2.25)\psline{-}(4.75,2.25)(5.0,2.0)
\psline{-}(5.0,2.0)(5.25,2.25) \psline{-}(5.25,2.25)(5.5,2.5)
%\psline{-}(11.0,5.0)(11.5,4.5) \psline{-}(11.5,4.5)(12.0,4.0)
%\psline{-}(12.0,4.0)(12.5,3.5)% \psline{-}(13.5,2.0)(14.0,1.5)
%\psline{-}(14.0,1.5)(14.5,1.0) \psline{-}(14.5,1.0)(15.0,0.5)

\psline{-}(0.25,1.25)(0.5,1.0) \psline{-}(0.5,1.0)(0.75,0.75)
\psline{-}(0.75,0.75)(1.0,1.0) \psline{-}(1.0,1.0)(1.25,1.25)
 \psline{-}(1.25,1.25)(1.5,1.0)\psline{-}(1.5,1.0)(1.75,0.75)
\psline{-}(1.75,0.75)(2.0,0.5) \psline{-}(2.0,0.5)(2.25,0.25)
%\psline{-}(4.5,2.0)(5.0,1.5) \psline{-}(5.0,1.5)(5.5,1.0)
%\psline{-}(5.5,1.0)(6.0,0.5)% \psline{-}(13.5,2.0)(14.0,1.5)
%\psline{-}(14.0,1.5)(14.5,1.0) \psline{-}(14.5,1.0)(15.0,0.5)

%\rput(2.5,0.25){{\small $i$}} %\rput(4.5,0.25){{\small $i+6$}}
%\rput(9.0,0.25){{\small $i+2$}} %\rput(12.0,0.25){{\small $i+8$}}
\rput(1.25,0.0){{\small $i$}} \rput(5.0,0.0){{\small $i+2$}}
\psline{-}(3.5,1.25)(3.75,1.0) \psline{-}(3.75,1.0)(4.0,0.75)
\psline{-}(4.0,0.75)(4.25,0.5) \psline{-}(4.25,0.5)(4.5,0.25)
 \psline{-}(4.5,0.25)(4.75,0.5)\psline{-}(4.75,0.5)(5.0,0.75)
\psline{-}(5.0,0.75)(5.25,0.5) \psline{-}(5.25,0.5)(5.5,0.25)
%\psline{-}(11.0,0.5)(11.5,1.0) \psline{-}(11.5,1.0)(12.0,1.5)
%\psline{-}(12.0,1.5)(12.5,1.0) \psline{-}(12.5,1.0)(13.0,0.5)
%\psline{-}(14.0,1.5)(14.5,1.0) \psline{-}(14.5,1.0)(15.0,0.5)

\end{pspicture}
\end{center}
\end{figure}

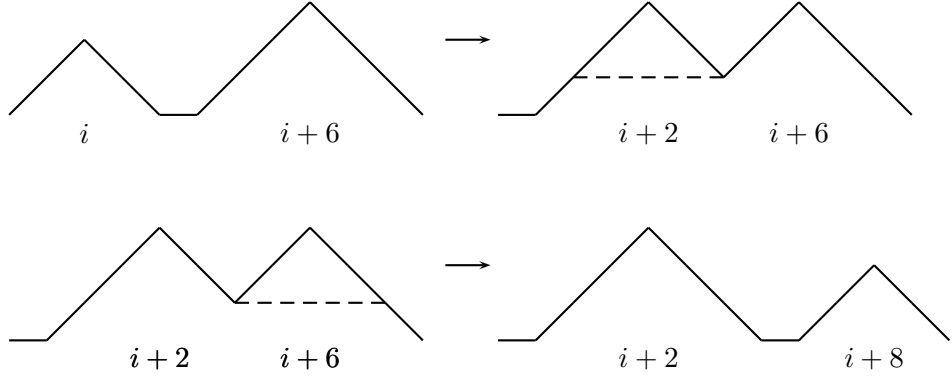
\begin{figure}[ht]
\caption{{\footnotesize Illustration of a sequence of moves where at some step there are two peaks of the same height in contact. The figure describes the motion of a peak of relative height 2 ``through''  a peak of relative height 3. The third step is simply a flip in the
identity of the two peaks.  Since both have the same height, the
second one is moved in the last step.}} \label{figure5}
\begin{center}
\begin{pspicture}(0,0)(13.0,5.5)
%axis
%\psline{->}(0.5,0.5)(0.5,2.0) \psline{->}(0.5,0.5)(13.0,0.5)
%units

%\psline{-}(12.5,0.5)(12.5,0.6) \psline{-}(13.0,0.5)(13.0,0.6)
%\psline{-}(13.5,0.5)(13.5,0.6) \psline{-}(14.0,0.5)(14.0,0.6)
%\psline{-}(14.0,0.5)(14.0,0.6)
%\psline{-}(14.5,0.5)(14.5,0.6)
\rput(1.5,3.25){{\small $i$}} \rput(4.5,3.25){{\small $i+6$}}
\rput(9.0,3.25){{\small $i+2$}} \rput(11.0,3.25){{\small $i+6$}}
%\rput(2.5,0.25){{\small $2$}}\rput(9.5,0.25){{\small$9$}}
%\rput(11.5,0.25){{\small $11$}}
 %\psline{-}(0.5,1.0)(0.6,1.0)
%\psline{-}(0.5,1.5)(0.6,1.5) %\psline{-}(0.5,2.0)(0.6,2.0)
%\psline{-}(0.5,2.5)(0.6,2.5) %\psline{-}(0.5,3.0)(0.6,3.0)
%\rput(0.25,1.5){{\small $1$}} %\rput(0.25,2.5){{\small $2 $}}
%graphic
\psline{-}(0.5,3.5)(1.0,4.0) \psline{-}(1.0,4.0)(1.5,4.5)
\psline{-}(1.5,4.5)(2.0,4.0) \psline{-}(2.0,4.0)(2.5,3.5)
\psline{-}(2.5,3.5)(3.0,3.5) \psline{-}(3.0,3.5)(3.5,4.0)
\psline{-}(3.5,4.0)(4.0,4.5) \psline{-}(4.0,4.5)(4.5,5.0)
\psline{-}(4.5,5.0)(5.0,4.5) \psline{-}(5.0,4.5)(5.5,4.0)
\psline{-}(5.5,4.0)(6.0,3.5) \psline{->}(6.3,4.5)(6.9,4.5)

\psline[linestyle=dashed]{-}(8.0,4.0)(10.0,4.0)
\psline[linestyle=dashed]{-}(3.5,1.0)(5.5,1.0)
\psline{->}(6.3,1.5)(6.9,1.5)

\rput(2.5,0.25){{\small $i+2$}} \rput(4.5,0.25){{\small $i+6$}}
%\rput(9.0,3.25){{\small $i+1$}} \rput(11.0,3.25){{\small $i+3$}}
\psline{-}(7.0,3.5)(7.5,3.5) \psline{-}(7.5,3.5)(8.0,4.0)
\psline{-}(8.0,4.0)(8.5,4.5) \psline{-}(8.5,4.5)(9.0,5.0)
 \psline{-}(9.0,5.0)(9.5,4.5)\psline{-}(9.5,4.5)(10.0,4.0)
\psline{-}(10.0,4.0)(10.5,4.5) \psline{-}(10.5,4.5)(11.0,5.0)
\psline{-}(11.0,5.0)(11.5,4.5) \psline{-}(11.5,4.5)(12.0,4.0)
\psline{-}(12.0,4.0)(12.5,3.5)% \psline{-}(13.5,2.0)(14.0,1.5)
%\psline{-}(14.0,1.5)(14.5,1.0) \psline{-}(14.5,1.0)(15.0,0.5)

\psline{-}(0.5,0.5)(1.0,0.5) \psline{-}(1.0,0.5)(1.5,1.0)
\psline{-}(1.5,1.0)(2.0,1.5) \psline{-}(2.0,1.5)(2.5,2.0)
 \psline{-}(2.5,2.0)(3.0,1.5)\psline{-}(3.0,1.5)(3.5,1.0)
\psline{-}(3.5,1.0)(4.0,1.5) \psline{-}(4.0,1.5)(4.5,2.0)
\psline{-}(4.5,2.0)(5.0,1.5) \psline{-}(5.0,1.5)(5.5,1.0)
\psline{-}(5.5,1.0)(6.0,0.5)% \psline{-}(13.5,2.0)(14.0,1.5)
%\psline{-}(14.0,1.5)(14.5,1.0) \psline{-}(14.5,1.0)(15.0,0.5)

\rput(2.5,0.25){{\small $i+2$}} \rput(4.5,0.25){{\small $i+6$}}
\rput(9.0,0.25){{\small $i+2$}} \rput(12.0,0.25){{\small $i+8$}}
\psline{-}(7.0,0.5)(7.5,0.5) \psline{-}(7.5,0.5)(8.0,1.0)
\psline{-}(8.0,1.0)(8.5,1.5) \psline{-}(8.5,1.5)(9.0,2.0)
 \psline{-}(9.0,2.0)(9.5,1.5)\psline{-}(9.5,1.5)(10.0,1.0)
\psline{-}(10.0,1.0)(10.5,0.5) \psline{-}(10.5,0.5)(11.0,0.5)
\psline{-}(11.0,0.5)(11.5,1.0) \psline{-}(11.5,1.0)(12.0,1.5)
\psline{-}(12.0,1.5)(12.5,1.0) \psline{-}(12.5,1.0)(13.0,0.5)
%\psline{-}(14.0,1.5)(14.5,1.0) \psline{-}(14.5,1.0)(15.0,0.5)

\end{pspicture}
\end{center}
\end{figure}

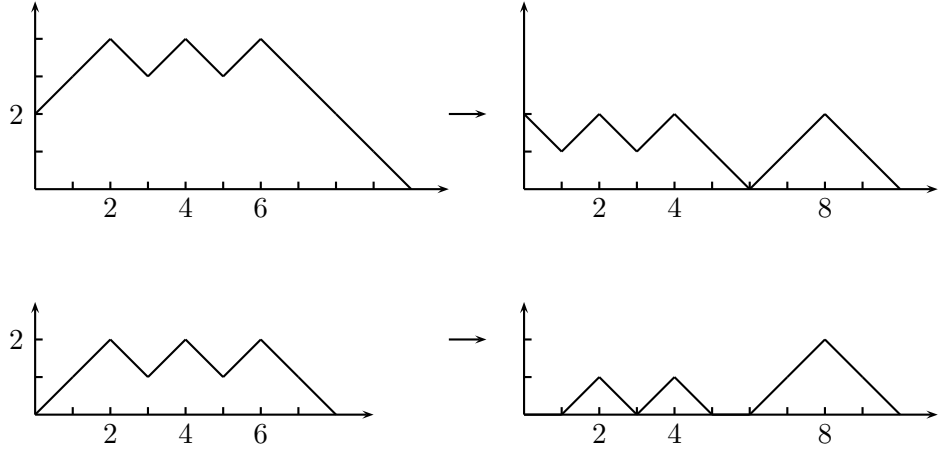
\begin{figure}[ht]
\caption{{\footnotesize Illustration of the displacement of the
first peak of relative height $2$ followed by two peaks of relative height $1$ with
the initial vertical  position being $2$ and $0$ respectively.}} \label{figure6}
\begin{center}
\begin{pspicture}(0,0)(13.0,7.5)
%axis
\psline{->}(0.5,3.5)(0.5,6.0) \psline{->}(0.5,3.5)(6.0,3.5)
\psline{->}(7.0,3.5)(7.0,6.0) \psline{->}(7.0,3.5)(12.5,3.5)
%units

%lines
\psline{-}(7.5,3.5)(7.5,3.6) \psline{-}(8.0,3.5)(8.0,3.6)
\psline{-}(8.5,3.5)(8.5,3.6) \psline{-}(9.0,3.5)(9.0,3.6)
\psline{-}(9.5,3.5)(9.5,3.6) \psline{-}(10.0,3.5)(10.0,3.6)
\psline{-}(10.5,3.5)(10.5,3.6) \psline{-}(11.0,3.5)(11.0,3.6)
\psline{-}(11.5,3.5)(11.5,3.6) %\psline{-}(12.0,0.5)(12.0,0.6)

\psline{-}(1.0,0.5)(1.0,0.6) \psline{-}(1.5,0.5)(1.5,0.6)
\psline{-}(2.0,0.5)(2.0,0.6) \psline{-}(2.5,0.5)(2.5,0.6)
\psline{-}(3.0,0.5)(3.0,0.6) \psline{-}(3.5,0.5)(3.5,0.6)
\psline{-}(4.0,0.5)(4.0,0.6) %\psline{-}(4.5,0.5)(4.5,0.6)

\psline{-}(1.0,3.5)(1.0,3.6) \psline{-}(1.5,3.5)(1.5,3.6)
\psline{-}(2.0,3.5)(2.0,3.6) \psline{-}(2.5,3.5)(2.5,3.6)
\psline{-}(3.0,3.5)(3.0,3.6) \psline{-}(3.5,3.5)(3.5,3.6)
\psline{-}(4.0,3.5)(4.0,3.6) \psline{-}(4.5,3.5)(4.5,3.6)
\psline{-}(5.0,3.5)(5.0,3.6)

\psline{-}(0.5,1.0)(0.6,1.0) \psline{-}(0.5,1.5)(0.6,1.5)
\psline{-}(0.5,4.0)(0.6,4.0) \psline{-}(0.5,4.5)(0.6,4.5)
\psline{-}(0.5,5.0)(0.6,5.0) \psline{-}(0.5,5.5)(0.6,5.5)
\psline{-}(7.0,1.0)(7.1,1.0) \psline{-}(7.0,1.5)(7.1,1.5)
\psline{-}(7.0,4.0)(7.1,4.0) \psline{-}(7.0,4.5)(7.1,4.5)
%\psline{-}(7.0,5.0)(7.1,5.0) \psline{-}(7.0,5.5)(7.1,5.5)
%\psline{-}(7.0,6.0)(7.1,6.0) \psline{-}(7.0,6.5)(7.1,6.5)
%\psline{-}(7.0,1.0)(7.1,1.0) \psline{-}(7.0,1.5)(7.1,1.5)

\rput(1.5,3.25){{\small $2$}} \rput(2.5,3.25){{\small $4$}}
\rput(3.5,3.25){{\small $6$}} \rput(8.0,3.25){{\small $2$}}
\rput(9.0,3.25){{\small $4$}} \rput(11.0,3.25){{\small $8$}}

\rput(0.25,4.5){{\small $2$}}

\rput(0.25,1.5){{\small $2$}}%\rput(9.5,0.25){{\small$9$}}
%\rput(11.5,0.25){{\small $11$}}
 %\psline{-}(0.5,1.0)(0.6,1.0)
%\psline{-}(0.5,1.5)(0.6,1.5) %\psline{-}(0.5,2.0)(0.6,2.0)
%\psline{-}(0.5,2.5)(0.6,2.5) %\psline{-}(0.5,3.0)(0.6,3.0)
%\rput(0.25,1.5){{\small $1$}} %\rput(0.25,2.5){{\small $2 $}}
\psline{-}(7.5,0.5)(7.5,0.6) \psline{-}(8.0,0.5)(8.0,0.6)
\psline{-}(8.5,0.5)(8.5,0.6) \psline{-}(9.0,0.5)(9.0,0.6)
\psline{-}(9.5,0.5)(9.5,0.6) \psline{-}(10.0,0.5)(10.0,0.6)
\psline{-}(10.5,0.5)(10.5,0.6) \psline{-}(11.0,0.5)(11.0,0.6)
\psline{-}(11.5,0.5)(11.5,0.6) %\psline{-}(12.0,0.5)(12.0,0.6)
%graphic
\psline{-}(0.5,4.5)(1.0,5.0) \psline{-}(1.0,5.0)(1.5,5.5)
\psline{-}(1.5,5.5)(2.0,5.0) \psline{-}(2.0,5.0)(2.5,5.5)
\psline{-}(2.5,5.5)(3.0,5.0) \psline{-}(3.0,5.0)(3.5,5.5)
\psline{-}(3.5,5.5)(4.0,5.0) \psline{-}(4.0,5.0)(4.5,4.5)
\psline{-}(4.5,4.5)(5.0,4.0) \psline{-}(5.0,4.0)(5.5,3.5)
%\psline{-}(5.5,4.0)(6.0,3.5) \psline{->}(6.3,4.5)(6.9,4.5)

%\psline[linestyle=dashed]{-}(8.0,4.0)(10.0,4.0)
%\psline[linestyle=dashed]{-}(3.5,1.0)(5.5,1.0)
\psline{->}(6.0,4.5)(6.5,4.5) \psline{->}(6.0,1.5)(6.5,1.5)

\rput(1.5,0.25){{\small $2$}} \rput(2.5,0.25){{\small $4$}}
\rput(3.5,0.25){{\small $6$}} %\rput(11.0,3.25){{\small $i+3$}}
\psline{-}(7.0,4.5)(7.5,4.0) \psline{-}(7.5,4.0)(8.0,4.5)
\psline{-}(8.0,4.5)(8.5,4.0) \psline{-}(8.5,4.0)(9.0,4.5)
 \psline{-}(9.0,4.5)(9.5,4.0)\psline{-}(9.5,4.0)(10.0,3.5)
\psline{-}(10.0,3.5)(10.5,4.0) \psline{-}(10.5,4.0)(11.0,4.5)
\psline{-}(11.0,4.5)(11.5,4.0) \psline{-}(11.5,4.0)(12.0,3.5)
%\psline{-}(12.0,4.0)(12.5,3.5)% \psline{-}(13.5,2.0)(14.0,1.5)
%\psline{-}(14.0,1.5)(14.5,1.0) \psline{-}(14.5,1.0)(15.0,0.5)

\psline{->}(0.5,0.5)(0.5,2.0) \psline{->}(0.5,0.5)(5.0,0.5)
\psline{->}(7.0,0.5)(7.0,2.0) \psline{->}(7.0,0.5)(12.5,0.5)

\psline{-}(0.5,0.5)(1.0,1.0) \psline{-}(1.0,1.0)(1.5,1.5)
\psline{-}(1.5,1.5)(2.0,1.0) \psline{-}(2.0,1.0)(2.5,1.5)
 \psline{-}(2.5,1.5)(3.0,1.0)\psline{-}(3.0,1.0)(3.5,1.5)
\psline{-}(3.5,1.5)(4.0,1.0) \psline{-}(4.0,1.0)(4.5,0.5)
%\psline{-}(4.5,2.0)(5.0,1.5) \psline{-}(5.0,1.5)(5.5,1.0)
%\psline{-}(5.5,1.0)(6.0,0.5)% \psline{-}(13.5,2.0)(14.0,1.5)
%\psline{-}(14.0,1.5)(14.5,1.0) \psline{-}(14.5,1.0)(15.0,0.5)

\rput(8.0,0.25){{\small $2$}} \rput(9.0,0.25){{\small $4$}}
\rput(11.0,0.25){{\small $8$}} %\rput(12.0,0.25){{\small $i+8$}}

%graph
\psline{-}(7.0,0.5)(7.5,0.5) \psline{-}(7.5,0.5)(8.0,1.0)
\psline{-}(8.0,1.0)(8.5,0.5) \psline{-}(8.5,0.5)(9.0,1.0)
 \psline{-}(9.0,1.0)(9.5,0.5)\psline{-}(9.5,0.5)(10.0,0.5)
\psline{-}(10.0,0.5)(10.5,1.0) \psline{-}(10.5,1.0)(11.0,1.5)
\psline{-}(11.0,1.5)(11.5,1.0) \psline{-}(11.5,1.0)(12.0,0.5)
%\psline{-}(12.0,1.5)(12.5,1.0) \psline{-}(12.5,1.0)(13.0,0.5)
%\psline{-}(14.0,1.5)(14.5,1.0) \psline{-}(14.5,1.0)(15.0,0.5)

\end{pspicture}
\end{center}
\end{figure}

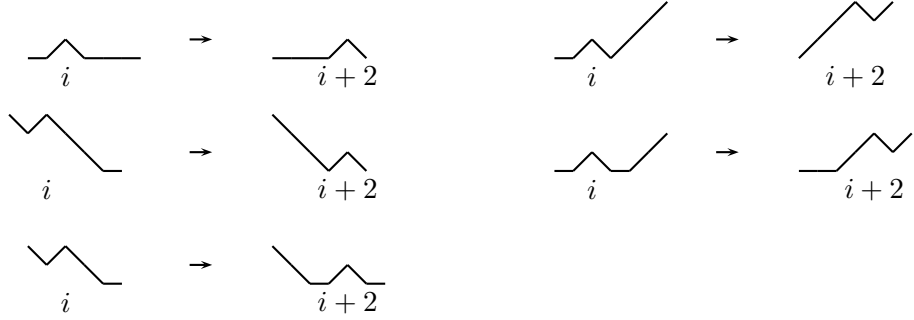
\begin{figure}[ht]
\caption{{\footnotesize Possible moves for a peak of relative height 1 within a path with all peaks being of relative height larger than 1.}}
\label{figure7}
\begin{center}
\begin{pspicture}(0,0)(12.0,5.0)
%axis
%\psline{->}(0.5,0.5)(0.5,2.0) \psline{->}(0.5,0.5)(13.0,0.5)
%units

% top right

\rput(7.75,3.0){{\small $i$}}
 \psline{-}(7.25,3.25)(7.5,3.25)\psline{-}(7.5,3.25)(7.75,3.5)
 \psline{-}(7.75,3.5)(8.0,3.25)\psline{-}(8.0,3.25)(8.25,3.5)
 \psline{-}(8.25,3.5)(8.5,3.75)\psline{-}(8.5,3.75)(8.75,4.0)
%\psline{-}(1.75,1.5)(2.0,1.0) \psline{-}(2.0,1.0)(2.25,0.5)
%\psline{-}(4.5,2.0)(5.0,1.5) \psline{-}(5.0,1.5)(5.5,1.0)
%\psline{-}(5.5,1.0)(6.0,0.5)% \psline{-}(13.5,2.0)(14.0,1.5)
%\psline{-}(14.0,1.5)(14.5,1.0) \psline{-}(14.5,1.0)(15.0,0.5)

\psline{->}(9.4,3.5)(9.7,3.5)% \psline{->}(2.9,2.5)(6.1,2.5)

%\rput(2.5,0.25){{\small $i$}} %\rput(4.5,0.25){{\small $i+6$}}
%\rput(9.0,0.25){{\small $i+2$}} %\rput(12.0,0.25){{\small $i+8$}}
 \rput(11.25,3.0){{\small $i+2$}} \psline{-}(10.5,3.25)(10.75,3.5)
\psline{-}(10.75,3.5)(11.0,3.75) \psline{-}(11.0,3.75)(11.25,4.0)
\psline{-}(11.25,4.0)(11.5,3.75)
 \psline{-}(11.5,3.75)(11.75,4.0)% \psline{-}(11.75,3.25)(12.0,3.25)
%\psline{-}(5.0,1.5)(5.25,1.0) \psline{-}(5.25,1.0)(5.5,0.5)
%\psline{-}(11.0,0.5)(11.5,1.0) \psline{-}(11.5,1.0)(12.0,1.5)
%\psline{-}(12.0,1.5)(12.5,1.0) \psline{-}(12.5,1.0)(13.0,0.5)
%\psline{-}(14.0,1.5)(14.5,1.0) \psline{-}(14.5,1.0)(15.0,0.5)

% middle right

\rput(7.75,1.5){{\small $i$}}
 \psline{-}(7.25,1.75)(7.5,1.75)\psline{-}(7.5,1.75)(7.75,2.0)
 \psline{-}(7.75,2.0)(8.0,1.75)\psline{-}(8.0,1.75)(8.25,1.75)
 \psline{-}(8.25,1.75)(8.5,2.0) \psline{-}(8.5,2.0)(8.75,2.25)
%\psline{-}(1.75,1.5)(2.0,1.0) \psline{-}(2.0,1.0)(2.25,0.5)
%\psline{-}(4.5,2.0)(5.0,1.5) \psline{-}(5.0,1.5)(5.5,1.0)
%\psline{-}(5.5,1.0)(6.0,0.5)% \psline{-}(13.5,2.0)(14.0,1.5)
%\psline{-}(14.0,1.5)(14.5,1.0) \psline{-}(14.5,1.0)(15.0,0.5)

\psline{->}(9.4,2.0)(9.7,2.0)% \psline{->}(2.9,2.5)(6.1,2.5)

%\rput(2.5,0.25){{\small $i$}} %\rput(4.5,0.25){{\small $i+6$}}
%\rput(9.0,0.25){{\small $i+2$}} %\rput(12.0,0.25){{\small $i+8$}}
 \rput(11.5,1.5){{\small $i+2$}}
\psline{-}(10.5,1.75)(10.75,1.75) \psline{-}(10.75,1.75)(11.0,1.75)
\psline{-}(11.0,1.75)(11.25,2.0) \psline{-}(11.25,2.0)(11.5,2.25)
 \psline{-}(11.5,2.25)(11.75,2.0) \psline{-}(11.75,2.0)(12.0,2.25)
%\psline{-}(5.0,1.5)(5.25,1.0) \psline{-}(5.25,1.0)(5.5,0.5)
%\psline{-}(11.0,0.5)(11.5,1.0) \psline{-}(11.5,1.0)(12.0,1.5)
%\psline{-}(12.0,1.5)(12.5,1.0) \psline{-}(12.5,1.0)(13.0,0.5)
%\psline{-}(14.0,1.5)(14.5,1.0) \psline{-}(14.5,1.0)(15.0,0.5)

% top left

\rput(0.75,3.0){{\small $i$}}
 \psline{-}(0.25,3.25)(0.5,3.25)\psline{-}(0.5,3.25)(0.75,3.5)
 \psline{-}(0.75,3.5)(1.0,3.25)\psline{-}(1.0,3.25)(1.25,3.25)
 \psline{-}(1.25,3.25)(1.5,3.25)\psline{-}(1.5,3.25)(1.75,3.25)
%\psline{-}(1.75,1.5)(2.0,1.0) \psline{-}(2.0,1.0)(2.25,0.5)
%\psline{-}(4.5,2.0)(5.0,1.5) \psline{-}(5.0,1.5)(5.5,1.0)
%\psline{-}(5.5,1.0)(6.0,0.5)% \psline{-}(13.5,2.0)(14.0,1.5)
%\psline{-}(14.0,1.5)(14.5,1.0) \psline{-}(14.5,1.0)(15.0,0.5)

\psline{->}(2.4,3.5)(2.7,3.5)% \psline{->}(2.9,2.5)(6.1,2.5)

%\rput(2.5,0.25){{\small $i$}} %\rput(4.5,0.25){{\small $i+6$}}
%\rput(9.0,0.25){{\small $i+2$}} %\rput(12.0,0.25){{\small $i+8$}}
 \rput(4.5,3.0){{\small $i+2$}}
\psline{-}(3.5,3.25)(3.75,3.25) \psline{-}(3.75,3.25)(4.0,3.25)
\psline{-}(4.0,3.25)(4.25,3.25) \psline{-}(4.25,3.25)(4.5,3.5)
 \psline{-}(4.5,3.5)(4.75,3.25) %\psline{-}(4.75,3.25)(5.0,3.25)
%\psline{-}(5.0,1.5)(5.25,1.0) \psline{-}(5.25,1.0)(5.5,0.5)
%\psline{-}(11.0,0.5)(11.5,1.0) \psline{-}(11.5,1.0)(12.0,1.5)
%\psline{-}(12.0,1.5)(12.5,1.0) \psline{-}(12.5,1.0)(13.0,0.5)
%\psline{-}(14.0,1.5)(14.5,1.0) \psline{-}(14.5,1.0)(15.0,0.5)

% middle left

\rput(0.5,1.5){{\small $i$}}
 \psline{-}(0.0,2.5)(0.25,2.25)
 \psline{-}(0.25,2.25)(0.5,2.5)\psline{-}(0.5,2.5)(0.75,2.25)
 \psline{-}(0.75,2.25)(1.0,2.0)\psline{-}(1.0,2.0)(1.25,1.75)
 \psline{-}(1.25,1.75)(1.5,1.75)%\psline{-}(1.5,2.0)(1.75,1.5)
%\psline{-}(1.75,1.5)(2.0,1.0) \psline{-}(2.0,1.0)(2.25,0.5)
%\psline{-}(4.5,2.0)(5.0,1.5) \psline{-}(5.0,1.5)(5.5,1.0)
%\psline{-}(5.5,1.0)(6.0,0.5)% \psline{-}(13.5,2.0)(14.0,1.5)
%\psline{-}(14.0,1.5)(14.5,1.0) \psline{-}(14.5,1.0)(15.0,0.5)

\psline{->}(2.4,2.0)(2.7,2.0)% \psline{->}(2.9,2.5)(6.1,2.5)

%\rput(2.5,0.25){{\small $i$}} %\rput(4.5,0.25){{\small $i+6$}}
%\rput(9.0,0.25){{\small $i+2$}} %\rput(12.0,0.25){{\small $i+8$}}
 \rput(4.5,1.5){{\small $i+2$}}
\psline{-}(3.5,2.5)(3.75,2.25) \psline{-}(3.75,2.25)(4.0,2.0)
\psline{-}(4.0,2.0)(4.25,1.75) \psline{-}(4.25,1.75)(4.5,2.0)
 \psline{-}(4.5,2.0)(4.75,1.75) %\psline{-}(4.75,1.75)(5.0,1.5)
%\psline{-}(5.0,1.5)(5.25,1.0) \psline{-}(5.25,1.0)(5.5,0.5)
%\psline{-}(11.0,0.5)(11.5,1.0) \psline{-}(11.5,1.0)(12.0,1.5)
%\psline{-}(12.0,1.5)(12.5,1.0) \psline{-}(12.5,1.0)(13.0,0.5)
%\psline{-}(14.0,1.5)(14.5,1.0) \psline{-}(14.5,1.0)(15.0,0.5)

% bottom left

\rput(0.75,0.0){{\small $i$}}
 \psline{-}(0.25,0.75)(0.5,0.5)\psline{-}(0.5,0.5)(0.75,0.75)
 \psline{-}(0.75,0.75)(1.0,0.5)\psline{-}(1.0,0.5)(1.25,0.25)
 \psline{-}(1.25,0.25)(1.5,0.25)%\psline{-}(1.5,2.0)(1.75,1.5)
%\psline{-}(1.75,1.5)(2.0,1.0) \psline{-}(2.0,1.0)(2.25,0.5)
%\psline{-}(4.5,2.0)(5.0,1.5) \psline{-}(5.0,1.5)(5.5,1.0)
%\psline{-}(5.5,1.0)(6.0,0.5)% \psline{-}(13.5,2.0)(14.0,1.5)
%\psline{-}(14.0,1.5)(14.5,1.0) \psline{-}(14.5,1.0)(15.0,0.5)

\psline{->}(2.4,0.5)(2.7,0.5)% \psline{->}(2.9,2.5)(6.1,2.5)

%\rput(2.5,0.25){{\small $i$}} %\rput(4.5,0.25){{\small $i+6$}}
%\rput(9.0,0.25){{\small $i+2$}} %\rput(12.0,0.25){{\small $i+8$}}
 \rput(4.5,0.0){{\small $i+2$}}
\psline{-}(3.5,0.75)(3.75,0.5) \psline{-}(3.75,0.5)(4.0,0.25)
\psline{-}(4.0,0.25)(4.25,0.25) \psline{-}(4.25,0.25)(4.5,0.5)
 \psline{-}(4.5,0.5)(4.75,0.25) \psline{-}(4.75,0.25)(5.0,0.25)
%\psline{-}(5.0,1.5)(5.25,1.0) \psline{-}(5.25,1.0)(5.5,0.5)
%\psline{-}(11.0,0.5)(11.5,1.0) \psline{-}(11.5,1.0)(12.0,1.5)
%\psline{-}(12.0,1.5)(12.5,1.0) \psline{-}(12.5,1.0)(13.0,0.5)
%\psline{-}(14.0,1.5)(14.5,1.0) \psline{-}(14.5,1.0)(15.0,0.5)

\end{pspicture}
\end{center}
\end{figure}

The generating function associated to these displacement is $[(q^2;q^2)_{m_{1}}(q^2;q^2)_{m_2}]^{-1}$. Collecting all factors, we recover the proper extension of the summand (\ref{hypo})  with the modes $m_1$ and $m_2$ inserted. With $n_j\rw {\tilde N}_j$, this is the expected result (cf. Proposition \ref{GFjp}).  By construction, when summed over all values of $m_j$, this function counts the number of paths that starts at $(0,2\ka+2-2i)$ with all  peaks at even $x$ coordinates and whose maximal  height is $\leq K+1$.
\end{proof}

 %=================================================================

 \section{Restricted ${\cal E}$-partitions}

 There is a simple transformation  relating  a jagged partition into a partition satisfying a non-increasing condition. One first doubles each part and replaces  every pair $(2r-2,2r)$ by the pair $(2r-1,2r-1)$.
 All entries are thereby ordered in non-increasing order from left to right and the length is preserved. By construction, every odd part must have even multiplicity. For instance
\begin{equation*}
(3,4,3,2,1,2,1,0,1) \quad \rw \quad \l(7,7,6,4,3,3,2,1,1\r)\;.
\end{equation*}
Let us call the resulting partition an ${\cal E}$-partition (where the ${\cal E}$ reminds of a build-in {\it eveness}).
The transformation of a jagged partition into an ${\cal E}$-partition is obviously a bijection. This bijection leads to the following equality.

\begin{proposition}\label{J=E} The numbers $J_{K,i}(n,m)$ and $ E_{K,i}(2n,m)$ defined in Theorems \ref{J=P} and \ref{burge} are equal.
\end{proposition}

\begin{proof} In the transformation from jagged partition to ${\cal E}$-partition, the weight is doubled,  the length is preserved and the frequency condition on the number of pairs of 01 (at most $i-1$) becomes a condition on the maximal number of 1, namely $2i-2$. The proof reduces then to show that for the ${\cal E}$-partition corresponding  to a given $K$-restricted jagged partition,  the restrictions (\ref{rjag}) become
\begin{equation}\label{simp}
p_j \geq  p_{j+K-1} +2 \;.
\end{equation}
%It is clear that $p_i\geq p_{i+K-1}$. One then only has to check that the following pairs of values of $(p_i,p_{i+K-1})$ are not allowed: $(2r,2r)$, $(2r+1,2r+1)$, $(2r+1,2r)$ and $(2r,2r-1)$. Suppose that $p_i= p_{i+K-1} = 2r$. This means that
%$ n_i=r,\; n_{i-1}\geq r$ and $n_{i+1}\leq r$.
% \qquad {\rm and}\qquad n_{i+K-1}= r
% \quad n_{i+K-2}\geq r, \quad n_{i+K}\leq r\; \end{equation}
%This possibility requires  $n_i = n_{i+K-1}$, in which case  we need to apply the  second  condition in  (\ref{rjag}). This forces $n_{i+1}= r+1$, which in turn violates the  inequality $n_{i+1}\leq r$. The other cases are treated similarly.
%\begin{equation}\label{rjag}
In the case where  $n_j \geq  n_{j+K-1} +1$, multiplying the parts by 2 and possibly rearranging them produces a difference $p_j -  p_{j+K-1}$ which is at least 2. If we have instead
$ n_j = n_{j+1}-1 =  n_{j+K-2}+1= n_{j+K-1} $, the transformation produces a difference which is precisely 2.
\end{proof}

%Note that there is a much more direct way of getting at (\ref{simp}). Instead of considering the defining conditions, one focuses  on the description of those partitions that are excluded. An example will illustrate the point. For $K=6$, the jagged partitions that  are excluded by (\ref{rjag}) are those which contain as a subsequence, a multiple of one of the following (omitting the comas between parts):
%\begin{equation*} (010101)\;, \;(110101)\;, \;(111101)\;, \;(111111)\;, \;(121111)\;, \;(121211).
%\end{equation*}
%The corresponding sub-$\cal E$-partitions are respectively multiple of
%\begin{equation*} (111111)\;, \;(221111)\;, \;(222211)\;, \;(222222)\;, \;(332222)\;, \;(333322),
%\end{equation*}
%and in all cases $p_i\leq p_{i+5}+1$.
%On the other hand, the two jagged subsequences $(120101)$ and $(121101)$ are allowed and these correspond to (331111) and (333311) respectively, which both  satisfy $p_i\geq p_{i+5}+2$.

The relation (\ref{simp})  points toward the  neat advantage of using ${\cal E}$-partitions:  not only these are genuine partitions, but the restrictions reduce to a simple `difference 2 condition
at distance $K-1$'.
In a sense, this enlightens  the somewhat mysterious  nature  of the original restrictions (\ref{rjag}).
%
% while the 3-restriced h-parttions of weight 3 are
% \begin{equation}
% \l\{(3),\; (2,1),\; \; \l(\frac32,\frac32\r) ,\;  \l(2,\frac12,\frac12\r)  \r\}
% \end{equation}

$\cal E$-partitions that are enumerated by $E_{K,i}(2n,m) $ have the following equivalent frequency  characterization. If $f_i$ stands for  the frequency of the part $i$ in the ${\cal E}$-partition, then $f_{2i-1}$ is even, $f_1\leq 2i-2$ and the restriction (\ref{simp}) becomes:
\begin{equation*}
f_i +f_{i+1}\leq  K-1.
\end{equation*}

% \begin{corollary}\label{GFE} The generaring fucntion
% for the restricted  $\cal E$ partition enumerated
% by  $ E_{K,i}(2n,m)$ defined in Theorem 1 is $G_{K,i}(z;q^2)$.
%  \end{corollary}

\let\a\alpha
\let\b\beta
%=================================================================

\section{Correspondence between paths and $\cal E$-partitions}

The equality $P_{K,i}(2n)= E_{K,i}(2n)$, that follows from Theorem \ref{J=P}  and Proposition \ref{J=E}, is proved  here by means of a bijection.
The bijection at work is precisely the one given by  Burge \cite{Bu} between a partition with a frequency condition and a two-word sequence, reinterpreted in terms of a path \cite{BreL}.  We then  only have to check that  $f_{2j-1}= 0 $ mod 2 if and only if the $x$-coordinates of the peaks are even.

The  Burge correspondence relies on the characterization of a partition in terms of non-overlapping pairs of adjacent frequencies $(f_j,f_{j+1})$ with $f_{j+1}>0$, starting the pairing from the largest part. For instance, for the following $\cal E$-partition of 62, we have the following pairing :
\begin{equation*}
(9,9,8,7,7,7,7,4,2,1,1) : \quad \begin{matrix}
i: \,0&\phantom{(}1&2\phantom{(}&\phantom{(}3
&4\phantom{(}&5&\phantom{(}6&7\phantom{(}&\phantom{(}8&9\phantom{(.} \\
  f_i: \, 0&(2&1)&(0&1)&0&(0&4)&(1&2). \end{matrix}
\end{equation*}
Let us now define a sequence of two operations on the set of paired frequencies. If $(f_0,f_1)$ is not a pair, we act with $\a$ defined as follows:
\begin{equation*}\a: (f_j,f_{j+1})\rw (f_j+1,f_{j+1}-1) \qquad \forall\, j\geq 1 .
\end{equation*}
If $(f_0,f_1)=(0,f_1)$ is a pair, we act with $\b$ defined as follows:
\begin{equation*}\b: \left\{ \begin{matrix}
&(0,f_1)\rw (0,f_1-1)\phantom{\qquad\qquad\;}& \\  & (f_j,f_{j+1})\rw (f_j+1,f_{j+1}-1) & \forall\, j>1 . \end{matrix} \right.
\end{equation*}
After each operation, the pairing is modified according to the new values of the frequencies. We then act successively with $\a$ or $\b$ on the partition until all frequencies become zero. The ordered sequence of $\a$ and $\b$ so obtained is then reinterpreted as a path starting at a prescribed initial position, by considering  $\a$ to be an horizontal  or a southeast step,  and $\b$  a northeast step \cite{BreL}. One then adds to the end of the sequence the number of $\a$ needed to reach the horizontal axis.
For the above example, the sequence is $\a\b\a\b^3\a^6\b^2\a^{21}\b^5$ and we add $\a^5$ at the end to make the corresponding path reach the $x$ axis assuming that it starts at the origin.  In terms of the data $\{(x_j;h_j ) \}$, which specify the position of the peaks $x_j$ and their relative height $h_j$, the path corresponding to this sequence is $\{(2;1)\, (6;3)\, (14;2)\,(40;5)\}$.
This correspondence between a sequence  or a path and an $\cal E$-partition is clearly invertible and provides the desired bijection.

The following reformulation of the  Burge correspondence, due to  Bressoud \cite{BreL} (Sect. 4), will readily establish the remaining point to be verified, namely the interrelation between the eveness of  both $x_j$ and $f_{2\ell-1}$.  To each pair $(x_j;h_j)$ characterizing a peak,  we associate the integers $s_j$ and $r_j$ via the equality
\begin{equation*}x_j = s_j h_j + r_j \;, \qquad 0\leq r_j< h_j\;.
\end{equation*}
$s_j$ and $s_j+1$ are then parts of a ${\cal E}$-partition with the following frequencies:
\begin{equation*}
f_{s_j}= h_j-r_j\;, \qquad  f_{s_j+1}=  r_j\;.
\end{equation*}
By construction, the sum of these parts is $x_j$ and $f_{s_j} + f_{s_j+1} = h_j \leq K-1$. It is also simple to check that the frequency of an odd part is even if $x_j$ is even. When $x_j$ is even, if $s_j$ is odd, then $h_j$ and $r_j$ have the same parity, so that the frequency of $s_j$ is even; on the other hand, if $s_j$ is even, $r_j$ is  even and it  is  the frequency of $s_j+1$. The proof the inverse statement is similar.

 With this correspondence, there is  a potential problem with the frequency condition if the value of $s_{j-1}$ is close to that of $s_j$. In such a case,  the sum of the  frequencies of two consecutive integers is no longer properly bounded.  When this is so, one has to make a shuffle \cite{BreL}, that is, to replace
 \begin{equation*}
  (x_{j-1};h_{j-1})\; (x_j;h_j)\rightarrow (x'_{j-1};h'_{j-1})\;  (x'_j;h'_j)\;,
    \end{equation*}
    where
 \begin{equation*}\label{shuf}
 x'_{j-1}= x_j-2h\;, \qquad   h'_{j-1} = h_j \;, \qquad x'_j= x_{j-1}+2h\;, \qquad h_j'= h_{j-1} \;,
   \end{equation*}
   and
 $h=  {\rm min}\;(h_{j-1},h_j)$.
The point here is only to observe that this  operation preserves the parity of the positions.

\begin{proposition}\label{PE} The numbers $P_{K,i}(2n)= \sum_m P_{K,i}(2n,m)$ and $ E_{K,i}(2n)$, defined respectively in Theorems \ref{J=P} and \ref{burge}, are equal. \end{proposition}

  %=================================================================

  \section{Correspondence between paths and partitions with prescribed successive ranks}

 %  there correspond a Ferrer diagram. Let  $d$ denote
 % the side of its  Durfee square  (the largest square
 % within the diagram) and denote the  Frobenius description of the partition as

  A partition $\la= (\la_1,\la_2,\cdots)$, whose conjugate is written $\la'$, has a Frobenius representation:
  \begin{equation*}
 \la=  \begin{pmatrix} s_1&s_2&\cdots &s_d\\ t_1&t_2&\cdots &t_d \end{pmatrix}\;,
    \end{equation*}
    with $s_i= \la_i-i$, $t_i= \la'_i-i $ and $d$ is the largest integer such that $\la_d\geq d$. Note that $s_1>s_2>\cdots >s_d$, $t_1>t_2>\cdots >t_d$ and $\sum_i\la_i = d+\sum_j(s_j+t_j)$.
The successive  ranks are defined as $SR(j)= s_j-t_j$.

The bijection between paths and partitions with prescribed successive ranks is directly lifted from \cite{BreL}. Let $a$ stands for the vertical position of the origin of the path and $o_j$ be the number of horizontal moves in the path at the left of the peak ($x_j,y_j)$. Then, to each peak we associate the pair of integers $(s_j,t_j)$ defined by
% \begin{equation*}
 \begin{align*}
 &o_j ~{\rm even:} &~ & s_j=\frac12(x_j-y_j+a) &~& t_j=\frac12(x_j+y_j-a-2)\cr
 &o_j ~{\rm odd:} &~& s_j=\frac12(x_j+y_j+a-1) &~& t_j=\frac12(x_j-y_j-a-1)\;.
  \end{align*}
 % \end{equation*}
  Since
  \begin{equation*}
  s_j+t_j+1=x_j
  \end{equation*}
  and $x_j$ is even, $SR(j)$  must be odd.
Further constraints follow from the  expressions for $y_j$ and the bounds $1\leq y_j\leq K-1$:
  % \begin{equation*}
 \begin{align*}
 &o_j ~{\rm even:} &~& 1\leq y_j=t_j-s_j+a+1\leq K-1 &~& \Rw \qquad SR(j)\in [a+2-K,a]\cr
 &o_j ~{\rm odd:} &~& 1\leq y_j=s_j-t_j-a  \leq K-1&~& \Rw  \qquad SR(j)\in [a+1,K-1+a]\;.
  \end{align*}
% \end{equation*}
These imply  that $SR(j)\in [a+2-K,a+K-1]$. With $a=2\ka-2i$, we have $SR(j)\in [2+\e-2i, 2K+\e-1-2i] $. Because $SR(j)$ must  be odd, this range can be shortened to $SR(j)\in [3-2i, 2K-1-2i] $.

 %This transformation is manifestly invertible.

\begin{proposition}\label{JR} The numbers $J_{K,i}(n)$ and $ R_{K,i}(2n)$, defined respectively in  Theorems \ref{J=P} and \ref{burge}, are equal. \end{proposition}

%=================================================================

\section{Related Rogers-Ramanujan identities}

Theorem 11 of \cite{FJM.E} and Proposition \ref{GFjp} imply the following Rogers-Ramanujan-type identity.

\begin{proposition} The multiple sum $G_{K,i}(1;q)$ defined in (\ref{defFabb}), with $K=2\kappa-\e$ and $1\leq i\leq
\kappa$  can be expressed in product form as
\begin{align*}\label{theor}
\sum_{m_{1},\dots,m_{\ka-1}=0}^\y  \frac{q^{\frac12({\tilde N}_{1}^2+\cdots+ {\tilde N}_{K-1}^2 +{\tilde M})+{\tilde L}_{2i}}
}{ (q)_{m_{1}}\cdots (q)_{m_{K-1}} }&= \prod_{n=1}^\infty  (1+ q^n)  \prod_{\substack{n\not
= 0,
\pm i\\
\;{\rm mod}\, (K +1)}}^\infty  \frac1{(1- q^n)}  & &   ( 2i<K+1 ) \cr
 &=
 \prod_{n\not = 0\;{\rm mod}\;
\kappa}^\infty \frac{(1+ q^n)}{ (1- q^n)}   && (\epsilon=1, i= \kappa) \;.
\end{align*}
\end{proposition}

The product form is manifestly the generating functions for those overpartitions enumerated by $O_{K,i}(n,m)$ defined in Theorem \ref{burge}. (With $q\rw q^2$ and using the identity $\prod (1+ q^n) = \prod (1-q^{2n-1})^{-1}$, this receives an alternative interpretation given in \cite{Bur}, Theorem 3.) This, together with Propositions \ref{J=E} and \ref{JR}, complete the proof of Theorem \ref{burge}.

%===================================================================

\section{Concluding remarks }

\subsection{Jagged partitions and  overpartitions}

There is a natural bijection between overpartitions and jagged
partitions, obtained as follows \cite{Lo}. Replace adjacent integers $(n,n+1)$ within
the jagged partition by $2n+1$ and similarly replace adjacent integers
$(n,n)$ by $2n$.  The remaining entries of the jagged partitions, necessarily distinct,
 are then overlined. (A similar bijection has also been obtained in \cite{Mah}.)

 Given this, it is natural to seek for  a possible relation between \cite{CM} and the present work. The authors of \cite{CM}  have considered a lattice path representation of overpartitions. The corresponding paths  are genuine generalizations of the usual paths of \cite{BreL}   in that some downward vertical moves are allowed. In contrast, our paths are standard ones, up to the restriction on the peak positions. Actually, the paths considered here can  more naturally be viewed as path-representations of $\cal E$-partitions rather than of jagged partitions per se.

The $k$-dependent restrictions considered in \cite{CM} are the following conditions `at distance $k$': $\la_i\geq \la_{i+k-1}+1 $ if the part $\la_{i+k-1}$ is overlined and $\la_i\geq \la_{i+k-1}+2 $ otherwise.  Are these  related to the restrictions (\ref{rjag})? Unfortunately not. The transformation between a jagged partition and an overpartition does not preserves the length, nor any notion of `distance' between parts.% which could have been shortened by a common factor say.

As pointed out in the introduction, the conditions (\ref{rjag}) arise naturally form a physical problem. It will be interesting to see  whether the above restricted overpartitions will also arise in a physical context.

\subsection{Generalized jagged partitions}

A natural axis of generalization is to consider other types of jagged partitions. The bijection between an unrestricted jagged partition and  an ordinary partition with generic frequency constraints (like an $\cal E$-partition) can be extended to  special generalized jagged partitions. For instance, we expect this to work for the   $0^p1$-partitions (where $p=0,1$ correspond respectively to ordinary partitions and the jagged partitions  studied here), defined by the conditions (cf. Definition 26 in Sect. 5.3 of \cite{FJM.R})
 \begin{equation*}
n_j\geq n_{j+s} -1\quad{\rm for}\quad  1\leq s\leq p, \qquad {\rm and} \quad n_j\geq n_{j+p+1} .
 \end{equation*}
 The procedure is to multiply all the parts of a $0^p1$-partition by $p+1$ and provide a rearrangement rule that ensures  the resulting parts to be  non-increasing.
 Take for instance $p=2$. The parts are multiplied by 3 and then rearranged as follows:
%   \begin{equation*}
      \begin{align*}
 (3r,3r,3r+3) &\rw (3r+1,3r+1,3r+1)\\ (3r,3r+3,3r+3) &\rw (3r+2,3r+2,3r+2)\\  (3r+3,3r,3r+3) & \rw (3r+3,3r+2,3r+1).
\end{align*}
% \end{equation*}
 A 001-jagged partition of weight $n$ is thus transformed into a partition of weight $3n$ satisfying the following frequency conditions (where $\e=0,1$):
   \begin{equation*}
f_{3j+1+\e} \equiv 0 \; {\rm mod}\, 3
\qquad {\rm or}\qquad
 f_{3j+1} \equiv  f_{3j+2} \equiv 1 \; {\rm mod}\, 3.
\end{equation*}
But the main problem with these extensions lies in the difficulty of finding restriction conditions  `at distance $k-1$'  that would  satisfy the  following naturalness criterion: the restriction excludes exactly one sequence of $k$ adjacent  parts for  each (allowed) value of the weight. No such condition has been found for $p>1$.

\vskip0.3cm
\noindent {\bf ACKNOWLEDGMENTS}

The work of PJ is supported by EPSRC and partially  by the EC
network EUCLID (contract number HPRN-CT-2002-00325), while that of  PM is supported  by NSERC.


\begin{thebibliography}{99}
\addcontentsline{toc}{section}{References}


% \bibitem{AB}
% A.K. Agarwal and D. Bressoud, {\it Latttice paths and
% hypergeometric series}. Pacific J. of Math. {\bf 136} (1989) 209-228.


\bibitem{AnB}
G.E. Andrews and D. Bressoud, {\it On  the Burge correspondence between partitions and binary words}. Rocky Mtn. J.  Math. Soc. {\bf 24} (1980) 225-233.



\bibitem{Andrr}
G.E. Andrews, {\it Multiple $q$-series}, Houston J. Math. {\bf 7} (1981) 11-22.



\bibitem{Andr}
G.E. Andrews, {\it The theory of
partitions}, Cambridge Univ. Press, Cambridge, UK, (1984).


\bibitem{BFJM}
L. B\'egin, J.-F. Fortin, P. Jacob and P. Mathieu,  {\it Fermionic characters
for graded parafermions}, Nucl. Phys {\bf B659} (2003) 365-386.

\bibitem{BP}
A. ÊBerkovich and P. Paule, {\it Lattice paths, $q$-multinomials and two variants of the Andrews-Gordon identities},  Ramanujan J. {\bf 5} (2002) 409--425.



\bibitem{BreL}
D. Bressoud, {\it Lattice paths and Rogers-Ramanujan identities}, in {\it Number Theory, Madras 1987},
ed. K. Alladi. Lecture Notes in Mathematics {\bf 1395} (1987) 140-172.


 \bibitem{Bu}
 W.H. Burge, {\it A correspondence between partitions related to generalizations of the Ramanujan-Rogers identities}, Discrete Math. {\bf 34} (1981) 9-15.

 \bibitem{Bur}
 W.H. Burge,  {\it A three-way correspondence between partitions}, Europ. J. Comb,. {\bf 3} (1982) 195-213.


\bibitem{CL}
S. Cortel and J. Lovejoy, {\it Overpartitions}, Trans. Amer. Math. Soc. 356 (2004), 1623-1635.

\bibitem{CM}
S. Cortel and O. Mallet, {\it Overpartitions, lattice paths and Rogers-Ramanujan identities}, math.CO/0601463.


\bibitem{FJM.R}
J.-F. Fortin, P. Jacob and P. Mathieu, {\it Jagged partitions}, Ramanujan J. {\bf 10} (2005) 215-235.


\bibitem{FJM.E}
J.-F. Fortin, P. Jacob and P. Mathieu, {\it Generating function for $K$-restricted jagged partitions}, Electronic J. Comb.
  {\bf 12} (2005) No 1, R12 (17 p.).


\bibitem{FJMa}
J.-F. Fortin, P. Jacob and P. Mathieu, {\it  SM(2,4$\kappa$)   fermionic characters
and restricted jagged partitions}, J. Phys. A: Math. Gen. {\bf 38} (2005) 1699-1709.






\bibitem{JM}
P. Jacob and P. Mathieu, {\it Graded parafermions: standard and quasi-particle bases}, Nucl. Phys.
{\bf B630} (2002) 433-452.



  \bibitem{JM.A}
  P. Jacob and P. Mathieu, {\it Parafermionic derivation of the  Andrews-type multiple sums} J. Phys. A: Math. Gen. {\bf 38} (2005) 8225-8238.

 % - 8238

  % {\it Parafermionic derivation of Andrews-type multiple sums}, hep-th/0505097

 \bibitem{Lo}
J. Lovejoy, {\it Constant terms, jagged partitions, and partitions with difference two at distance two}, to appear in Aequationes Mathematicae.




\bibitem{Mah}
K. Mahlburg, {\it The overpartition function modulo small powers of 2}, Discrete Math. {\bf 286} (2004), 263-267.

 \bibitem{Ole}
 S. O. Warnaar,
\textit{Fermionic solution of the Andrews--Baxter--Forrester model. I.
Unification of CTM and TBA methods},
J. Stat. Phys. \textbf{82} (1996), 657--685.




 \bibitem{War}
 S.O. Warnaar, {\it The generalized Borwein conjecture. II. Refined $q$-trinomial coefficients}, Discrete Math. {\bf 272} (2003) 215-258.

  \bibitem{Wa}
 S.O. Warnaar, private communication (2005).

  \end{thebibliography}
\end{document}